\documentclass[sn-mathphys-ay]{sn-jnl}

\usepackage{graphicx}%
\usepackage{multirow}%
\usepackage{amsmath,amssymb,amsfonts}%
\usepackage{amsthm}%
\usepackage{mathrsfs}%
\usepackage[title]{appendix}%
\usepackage{xcolor}%
\usepackage{textcomp}%
\usepackage{manyfoot}%
\usepackage{booktabs}%
\usepackage{algorithm}%
\usepackage{algorithmicx}%
\usepackage{algpseudocode}%
\usepackage{listings}%
\usepackage{soul}
\usepackage{ulem}
\usepackage{cancel}
\usepackage{comment}
\usepackage{array}
\usepackage{amsmath} 
\usepackage{subcaption}

\theoremstyle{thmstyleone}%
\newtheorem{theorem}{Theorem}
\theoremstyle{thmstyletwo}%

\theoremstyle{thmstylethree}%
\newtheorem{definition}{Definition}%

\newcommand{\ignore}[1]{}

\newcommand{\Rb}{\mathbb{R}}
\newcommand{\Eb}{\mathbb{E}}
\newcommand{\Zc}{\mathcal{Z}}
\newcommand{\Xb}{\mathbf{X}}
\newcommand{\xb}{\mathbf{x}}
\newcommand{\Yb}{\mathbf{Y}}
\newcommand{\yb}{\mathbf{y}}
\newcommand{\fsd}{\succeq_{(1)}}
\newcommand{\ssd}{\succeq_{(2)}}

\begin{document}

\title[Article Title]{Asset liability management under sequential stochastic dominance constraints}


\author[1]{\fnm{Giorgio} \sur{Consigli}}\email{giorgio.consigli@ku.ac.ae}

\author[2]{\fnm{Darinka} \sur{Dentcheva}}\email{darinka.dentcheva@stevens.edu}
\equalcont{These authors contributed equally to this work.}

\author*[3]{\fnm{Francesca} \sur{Maggioni}}\email{francesca.maggioni@unibg.it}
\equalcont{These authors contributed equally to this work.}

\author[3]{\fnm{Giovanni} \sur{Micheli}}\email{giovanni.micheli@unibg.it}
\equalcont{These authors contributed equally to this work.}

\affil*[1]{\orgdiv{Department of Mathematics}, \orgname{Khalifa University of Science and Technology}, \orgaddress{\street{Shakhbout Bin Sultan Street}, \city{Abu Dhabi}, \postcode{127788}, \country{UAE}}}

\affil[2]{\orgdiv{Department of Mathematics}, \orgname{Stevens Institute of Technology}, \orgaddress{\street{1 Castle Point Terrace}, \city{Hoboken}, \postcode{07030}, \state{New Jersey}, \country{USA}}}

\affil[3]{\orgdiv{Department of Management, Information and Production Engineering}, \orgname{University of Bergamo}, \orgaddress{\street{Via Einstein 2}, \city{Dalmine (BG)}, \postcode{24044}, \country{Italy}}}

\abstract{We consider a financial intermediary managing assets and liabilities exposed to several risk sources and seeking an optimal portfolio strategy to minimise the initial capital invested and the total risk associated with investment losses and financial debt.
We formulate the problem as a multistage stochastic programming model, with a time-consistent dynamic risk measure in the objective function to control the investment risk.
To ensure that the intermediary's financial equilibrium is preserved, we introduce a funding constraint in the model by enforcing in a time-consistent manner a sequential second-order stochastic dominance (SSD) of the portfolio return distribution over the liability distribution. 
We demonstrate that imposing the SSD constraint at the last-but-one stage is sufficient to enforce the SSD ordering at each stage.
To deal with the computational burden of associated MSP, we  develop a novel decomposition scheme integrating, for the first time in the literature,  time-consistent dynamic risk measures and sequential stochastic dominance constraints.
The proposed methodology is computationally validated on a case study developed on a property and casualty ALM problem.}



\maketitle

\section{Introduction}\label{sec1}

We consider a multi-period \textit{asset-liability management} (ALM) problem for a financial intermediary 
managing over a long-term horizon assets and liabilities exposed to several risk sources.
Inspired by a corporate case study, the ALM problem features a large insurance and financial intermediary, whose business structure is summarized by a technical division, an investment and a risk management division, which are responsible for \textit{liability} policies, \textit{asset} management and \textit{risk} assessment, respectively. 
The decision maker, based on the risk assessment and the current asset-liability portfolios, intends to determine a minimal capital endowment and an asset-liability strategy able to preserve the intermediary's financial equilibrium while hedging a complex set of risk sources.

Due to the dynamic nature of the problem and the high uncertainty involved, the ALM problem finds a natural formulation as a 
\textit{multistage stochastic program} (MSP).
Several contributions have been proposed in this area in recent years. 
To mention just a few: \cite{RePEc:gro:rugsom:00a52} address a pension fund problem considering the adjustments of the contribution rate of the sponsor and the reallocation of the investments in several asset classes at various stages, \cite{urban2004} analyse the risk capital allocation problem for an insurance portfolio, \cite{gatzert2008} combine fair pricing and capital requirements for non-life insurance companies, \cite{alessandri2010} develop the analysis on risk capital requirements for a banking intermediary, \cite{dhaene2012} establish optimal capital allocation principles in presence of several risk sources, \cite{maume2015risk} propose a risk management approach to capital allocation by financial intermediaries, \cite{DUARTE2017177} develop a MSP model for open pension schemes with a thorough representation of a risk-based regulation, \cite{consigli2018} analyse the implications of capital constraints on the optimization of risk-adjusted returns in a dynamic model, \cite{lauria2022optimal} address a defined-benefit pension fund problem combining a stochastic control approach with a chance constraint on the funding ratio.

An important feature of the ALM problem is the high level of risk, which needs to be properly accounted for.
Coherent risk measures are popular approaches used in the literature to control risk in sequential decision problems. 
We refer to \cite{shapiro2021lectures} and \cite{Follmer2008} for an overview on coherent risk measures and to \cite{chen2017multi} for a discussion on their application to the financial context.

Additionally, a very popular and powerful tool for risk assessment and decision-making in MSP models is provided by \textit{stochastic dominance} (SD), which helps in comparing different uncertain prospects or investment options by evaluating their risk and return characteristics.
Optimization problems with stochastic dominance constraints were first introduced in \cite{dentcheva2003G} and further developed in \cite{Dentcheva2004}. 
While a rich literature analyses stochastic order relations for scalar random variables and their implications for decision making under uncertainty, much fewer works are devoted to the comparison of sequences, although dynamical systems are prevalent in practice. 
Some recent works include \cite{dentcheva2008stochastic,dentcheva2022risk,Escudero2016,Escudero2017,haskell2013stochastic}.  
The challenge in such comparisons arise from time-consistency considerations, which are crucial in the context of sequential decision making. 
The majority of sequential comparisons proposed in the literature lead to time-inconsistent decision problems (see \cite{Ruszcz2010}). Roughly speaking, a policy is defined as time-consistent if and only if the future planned decisions are actually implemented.

We emphasize that the notion of \textit{sequential stochastic dominance} differs substantially from the notion of \textit{multivariate stochastic dominance}.
The latter extends comparisons from scalar outcomes to random vectors, see \cite{dentcheva2009optimization,dentcheva2015optimization,armbruster2015models,dentcheva2016two}. 
These works apply the multivariate dominance relations to the entire distributions of the respective random vectors or processes.
Unlike sequential comparison methods, multivariate stochastic dominance does not explicitly incorporate the sequential structure of information, and it does not neglect events that cannot happen in the future.
Instead, it evaluates the entire process at once, without filtering out scenarios that become irrelevant at the time of comparison. Consequently, while the approach in the cited works is well-suited for high-level strategic planning, it does not ensure time consistency, a crucial property for decision-making frameworks that require non-contradictory policy.

In the ALM setting, SD provides a way to compare different decision policies or strategies under uncertainty, evaluating whether one policy consistently offers better outcomes than another, considering both risk and return.
An early application of SD criteria in a multistage ALM problem was due to \cite{yang2010asset} with a focus on risk control at specific stages. More recently, \cite{consigli2020long} solved an individual ALM problem over a long-term horizon, enforcing SD constraints at individual stages, under an independence assumption. 
More precisely, the SD relation was used to compare the performance of the constructed portfolio to a benchmark portfolio stage-wise at selected stages.  
A similar approach is proposed in \cite{kopa2018individual}, \cite{moriggia2019} and \cite{RePEc:spr:comgts:v:15:y:2018:i:2:d:10.1007_s10287-018-0299-8}. 
However, 
the proposed ALM models are either static or use time-inconsistent comparisons in a multistage setting. 
More recently, \cite{kopa2023multistage} incorporate multivariate stochastic dominance constraints into a pension funding problem, ensuring that the optimal strategy stochastically dominates a benchmark portfolio. However, as in other multivariate approaches, 
the comparison is performed over the entire stochastic process without explicitly considering the sequential information structure.

From a computational perspective, the application of stochastic dominance and risk measures further increases the computational burden of multistage ALM problems, requiring the development of suitable numerical approaches to solve real-case instances. 
Numerical methods for static optimization problems involving stochastic dominance relations as constraints or as multivariate objectives are proposed in \cite{dentcheva2003G,noyan2008valid,rudolf2008optimization,roman2006portfolio,luedtke2008,dentcheva2010inverse,fabian2011processing}.
In \cite{gulten2015two} an extended two-stage problem formulation with coherent measures of risk is presented and a numerical method for solving the problem is proposed, while a scenario decomposition method for multi-stage stochastic programming problems with coherent measures of risk is introduced in \cite{collado2012scenario}.

In this paper, we take the perspective of a financial intermediary managing assets and liabilities exposed to several risk sources and we formulate a MSP to determine an optimal portfolio strategy, controlling the investment risk by means of a time-consistent dynamic risk measure in the objective.
In order to ensure that the intermediary's financial equilibrium is preserved, we introduce a funding constraint in the model by enforcing in a time-consistent manner a dynamic second order stochastic dominance (SSD) of the portfolio return distribution over the liability distribution. 
The time-consistent stochastic comparison we introduce in the problem is new and differs from the one presented in the literature (see \cite{dentcheva2022risk}) since our stochastic comparison is required only  to a one-step look-ahead without the need to benchmark the recourse function sequences. 
We further demonstrate that imposing the SSD constraint at the last-but-one stage is sufficient to enforce the SSD ordering at each stage.
By adopting the proposed model, the \textit{ALM manager} may face some scenarios where the funding condition may actually worsen, but overall she/he will preserve an effective liability hedge.

Given the computational burden of MSP with risk measures and stochastic dominance constraints, in this article we also develop a decomposition scheme to efficiently solve the problem integrating, for the first time in the literature, a suitable numerical treatment for risk measures and sequential stochastic dominance constraints. 
Specifically, we propose a new version of the multi-cut method  from \cite{gulten2015two} in which additional event cuts approximate the time consistent stochastic order constraints and further cuts approximate the risk measures in the objective functions.

To sum up, the contributions of our paper can be summarized as follows. 

\begin{itemize} 
\item A novel multi-stage ALM model is presented, including time-consistent stochastic-ordering relations and a dynamic measure of risk for constructing optimal investment strategies under several risk sources.
Specifically, the enforcement of a sufficient funding condition through second order stochastic dominance in a multistage and time-consistent manner has not been proposed before and represents a key contribution of this work.
\item The imposition of the SSD constraint at the last-but-one stage is proven to be sufficient to enforce the SSD ordering at every stage (see Theorem \ref{theorem1}).
\item A new decomposition method for multistage stochastic optimization problems is developed, which includes risk measures as well as stochastic dominance constraints.
\item An extended set of financial-based validation evidence and sensitivity results, finally, analyze the impact of the stochastic dominance conditions on investment strategies and funding status.
\end{itemize}

The paper is organized as follows. We introduce the notions and properties  related to  stochastic dominance in Section \ref{sec:TC-MSD}. The ALM problem is then formulated in Section \ref{sec:ALM}. In Section \ref{sec:statmodel}, we present the stochastic models adopted to derive the full set of coefficients of the ALM model. The numerical approach developed to solve the ALM problem is presented in Section \ref{sec:decomposition}. Section \ref{sec:computations} discusses the computational evidences collected to validate the methodology and to support the decision-making process. Finally, conclusions are drawn in Section \ref{Conclusions}.
We leave to the appendices the detailed description of the stochastic models supporting scenario generation and the algorithm developed.

\section{Sequential stochastic dominance} \label{sec:TC-MSD}

First, we introduce the notions of stochastic dominance of first and higher order. In the following, the right-continuous cumulative distribution function (CDF) $F_Z(\eta)$ of a random variable $Z$ is defined as $F_Z(\eta)=P(Z\leq\eta)$ and the survival function of $Z$ is given by $\bar{F}_Z(\eta)= P (Z>\eta).$
The integrated distribution function $F^{(2)}_Z(\eta)$ is defined as follows:
\begin{equation*}
F^{(2)}_Z(\eta)=\int_{-\infty}^\eta F_Z(t)\,dt \,\text{ for }\eta\in\Rb.
\end{equation*}
Clearly, the function $F^{(2)}_Z(\cdot)$ is finite everywhere whenever $Z$ is integrable and it is convex as an integral of a non-decreasing function.
For a random variables $Z$ with a finite $k$-th moment, $k\geq 2$, we define recursively the functions
\begin{equation}
\label{kth}
 F^{(k+1)}_Z(\eta) = \int_{-\infty}^{\eta} F^{(k)}_Z(\alpha) \ d\alpha
 \quad \mbox{for} \  \eta \in \Rb.
 \end{equation}
\begin{definition}
\begin{itemize}
    \item[(i)] A random variable $V$ is stochastically larger than a random variable $Z$ with respect to the first order
stochastic dominance (denoted $V\fsd Z$) if $F_V(\eta)\leq F_Z(\eta)\,\text{ for all }\eta\in\Rb$.
\item[(ii)] For two random variables $V$ and $Z$, it is said that the variable $V$ is stochastically larger than $Z$  with respect to the $k$-th order stochastic dominance
(denoted $V \succeq_{(k)} Z$)  if
$F^{(k)}_V(\eta)\leq F^{(k)}_Z(\eta)$ for all $\eta\in\Rb$.
\end{itemize}
\end{definition}
Notice that the relation $V\fsd Z$ is also equivalent to $\bar{F}_Z(\eta)\leq \bar{F}_V(\eta)\,\text{ for all }\eta\in\Rb$,
meaning that $V$ takes larger values more frequently but comparisons of integrated survival functions lead to different relations than the $k$-order dominance.

We use the shorthand notation $a_+=\max(0,a)$ for any $a\in\Rb$. Changing the order of integration in the definition of the function $F^{(2)}_Z(\cdot)$, we obtain the following equivalent representation of the second-order relation: 
\begin{align}
V\ssd Z & \quad \Leftrightarrow \quad 
\Eb[\eta -V]_+\leq  \Eb[\eta -Z]_+,\quad \text{for all }\eta\in\Rb. \label{shortfall}
\end{align}
The second-order stochastic dominance relation can also be characterized by the respective quantile functions, which turned out to be very useful. 
Let $F^{-1}_Z(\cdot)$ be the left continuous inverse of the cumulative distribution function $F_Z(\cdot)$ defined by
\[
F^{-1}_Z(p)=\inf\{\eta: F_Z(\eta)\geq p\},\text{ for } 0<p<1.
\]
The absolute Lorenz function $L_Z:[0,1]\rightarrow \Rb$, introduced in the seminal work of Lorenz, see \cite{Lor1905},
is defined as the cumulative quantile function:
\begin{equation*}
\label{lfe}
L_Z(p)=\int_0^p F^{-1}_Z(t)dt \text{ for } 0<p\leq 1.
\end{equation*}
The definition of the function beyond the interval $(0,1]$ is extended by setting $L_Z(0)=0$ and
$L_Z(p)=\infty$ for $p\not\in [0,1]$. The Lorenz function is widely used in economics for comparison of income streams.

Interestingly the integrated distribution function and the Lorenz function are related via conjugate duality. It is shown in \cite{ogryczak2002} that
${L}_Z(\cdot)$ and $F^{(2)}_Z(\cdot)$ are Fenchel conjugate functions.
This result implies that relating the Lorenz functions of two integrable random variables provides an equivalent characterizations of the stochastic ordering relations, i.e.,
\begin{gather}
V\succeq_{(2)} Z\,\Leftrightarrow\, L_V(p)\geq L_Z(p)\quad \text{ for all } p\in[0,1] \label{lfso}.
\end{gather}
It is clear that $ Z\fsd V$ if and only if  $F^{(-1)}_Z(\eta)\geq F^{(-1)}_V(\eta)\,\text{ for all }\eta\in\Rb$ but a
quantile characterization for the relations of order $k>2$ is not available. \\


Let us turn to comparison of sequences. Given probability space $(\Omega,\mathcal{F},P)$, a filtration
$\mathcal{F}_1 \subset \dots \subset \mathcal{F}_{T}$,
with $\mathcal{F}_0=\{\emptyset,\Omega\}$ and $\mathcal{F}_{T}=\mathcal{F}$, denote
$ \Zc=\mathcal{L}_1(\Omega,\mathcal{F}_{1},P)\times\dots
\times \mathcal{L}_1(\Omega,\mathcal{F}_{T},P)$. We assume that the filtration is generated by the random data process $\{ \xi_t \}_{t=1}^T$ and denote the history of the data process until time $t$ by  $\xi_{[t]}$.

We wish to compare two sequences $\Xb=(X_1,X_2,\dots,X_{T})$ and $\Yb=(Y_1,Y_2,\dots,Y_{T})$ in $\Zc$ at any time $t=1,\dots T$ in a consistent way.
To this aim, we introduce the cumulative sum of $\Xb$ until time $t$ as a function of the history path $\xi_{[t]}$, i.e.,  
\[
\xb_t (\xi_{[t]}) = (X_1+X_2+\dots + X_t) (\xi_{[t]}),\quad t=1,\dots T.
\]
We shall denote the projected future value for the sequence $\Xb$ at time $t$ when 
$\xi_{[t]}$ is fixed defined as follows:
\[
\Xb_{t+1}| \xi_{[t]} =  X_{t+1}| \xi_{[t]} + \Eb_{t+1}\Big[ X_{t+2}|\xi_{[t]}+ 
\;\Eb_{t+2}\big[ X_{t+3}|\xi_{[t]} + \cdots + \Eb_{T-1}[X_T| \xi_{[t]}]\big]\Big].
\]
A time-consistent comparison, proposed in  \cite{dentcheva2022risk}, is obtained in the following way. 
We call the sequence $\Xb\in\Zc$ stochastically larger than the sequence $\Yb\in\Zc$ 
 if at any time $t$ and history $ \xi_{[t]}$ the following holds
 \begin{equation*}
 \xb_t(\xi_{[t]}) +\Xb_{t+1}\,| \xi_{[t]}\succeq_{\xi_{[t]}}\; \yb_t(\xi_{[t]}) + \Yb_{t+1}  \,| \xi_{[t]}, 
 \end{equation*}
 where the comparison $\succeq_{\xi_{[t]}}$ is an appropriately chosen stochastic ordering for scalar-valued random variables, which may be chosen depending on the available information at time $t$. The choice may depend on the current state of the data process, or simply on the time of comparison. In particular, we may postulate 
 \begin{equation*}
 \xb_t(\xi_{[t]}) +\Xb_{t+1}\,| \xi_{[t]}\ssd\; \yb_t(\xi_{[t]}) + \Yb_{t+1}  \,| \xi_{[t]}. 
 \end{equation*}

We denote the accumulated difference between $\Xb$ and $\Yb$ until time $t$ along the path (scenario) $\xi_{[t]}$ by $\sigma_t(\xi_{[t]})$, i.e.,
$\sigma_t(\xi_{[t]}) = \xb_t(\xi_{[t]}) - \yb_t(\xi_{[t]}).$
The relation in \eqref{long-ssd} is given by
\[
\int_0^p F^{-1}_{\xb_t(\xi_{[t]}) +\Xb_{t+1}\,| \xi_{[t]}}(t)dt
-\int_0^p F^{-1}_{\yb_t(\xi_{[t]}) +\Yb_{t+1}\,| \xi_{[t]}}(t)dt\geq 0
\]
for all $0<p\le 1$.
We transform the right-hand side as follows
\begin{eqnarray*}
& & \int_0^p F^{-1}_{\Xb_{t+1}\,| \xi_{[t]}}(t) + \xb_t(\xi_{[t]})\; dt -
\int_0^p F^{-1}_{\Yb_{t+1}\,| \xi_{[t]}}(t) +\yb_t(\xi_{[t]})\; dt\\
& = & \int_0^p F^{-1}_{\Xb_{t+1}\,| \xi_{[t]}}(t) - F^{-1}_{\Yb_{t+1}\,| \xi_{[t]}}(t) + \xb_t(\xi_{[t]}) - \yb_t(\xi_{[t]})\; dt \\
 & = & \int_0^p F^{-1}_{\Xb_{t+1}\,| \xi_{[t]}}(t) - F^{-1}_{\Yb_{t+1}\,| \xi_{[t]}}(t) + \sigma_t(\xi_{[t]})\; dt.
\end{eqnarray*}
Hence, \eqref{long-ssd} is equivalent to the following relation  at time $t$
 \begin{equation}
 \label{long-ssd}
 \sigma_t(\xi_{[t]}) +\Xb_{t+1}\,| \xi_{[t]}\ssd\;  \Yb_{t+1}  \,| \xi_{[t]}. 
 \end{equation}

This order satisfies the following definition of time-consistent sequential comparison.
\begin{definition} \label{def:def2}
A stochastic order between elements of  $\Zc$ is called time consistent, if for all sequences  $\Xb,\Yb \in \Zc$, for all  $1 \leq t < T-1$, and for all $\xi_{[t]}$ the relation $\Xb | (\xi_{[t]},\xi_{t+1}) \succeq_{\xi_{[t+1]}} \Yb | (\xi_{[t]},\xi_{t+1})$ for almost all\footnote{''Almost all'' is understood with respect to the conditional probability, given $\xi_{[t+1]}.$} $\xi_{t+1}$ implies that $\Xb | \xi_{[t]}\succeq_{\xi_{[t+1]}} \Yb | \xi_{[t]}$.  
\end{definition}

Definition \ref{def:def2}  is in harmony with the notion of time consistency of dynamic risk measures.
We shall show the following implication for a multi-stage stochastic programming problem with stochastic dominance relation of the second order imposed according to \eqref{long-ssd}.

\setcounter{theorem}{0}
\begin{theorem}
\label{theorem1}
The dynamic order \eqref{long-ssd} holds provided
\begin{equation}
\sigma_{T-1} (\xi_{[T-1]}) + \Xb_T | \xi_{[T-1]}\ssd \Yb_T | \xi_{[T-1]}    
\end{equation}
holds for almost all $\xi_{[T-1]}$. 
\end{theorem}

\begin{proof}
The assumption implies that for any history $\xi_{[T-1]} = (\xi_{[T-2]},\xi_{T-1})$:
 \begin{equation*}
 \xb_{T-2}(\xi_{[T-2]}) + X_{T-1}(\xi_{[T-1]}) + \Xb_{T}\,| \xi_{[T-1]} \ssd\; \yb_{T-2}(\xi_{[T-2]}) + Y_{T-1}(\xi_{[T-1]}) + \Yb_{T} | \, \xi_{[T-1]}. 
 \end{equation*}
 The second order stochastic dominance relation between two random variables implies that their expected values are related in the same way. Hence, for any history $\xi_{[T-1]}$
  \begin{equation*}
 \mathbb{E}_{T-1} [\xb_{t-2}(\xi_{[T-2]}) + ( X_{T-1} + \Xb_{T} ) | \xi_{[T-1]} ] \geq \mathbb{E}_{T-1} [\yb_{t-2}(\xi_{[T-2]}) + (Y_{T-1} + \Yb_{T}) | \xi_{[T-1]}]. 
 \end{equation*}
 Fixing a specific history $\xi_{[T-2]}$, we obtain
 \begin{equation*}
 \xb_{T-2}(\xi_{[T-2]}) + \mathbb{E}_{T-1} [ (X_{T-1} + \Xb_{T}) | \xi_{[T-2]} ] \geq \yb_{T-2}(\xi_{[T-2]}) + \mathbb{E}_{T-1} [(Y_{T-1} + \Yb_{T}) | \xi_{[T-2]}], 
 \end{equation*}
which is the same as
\begin{equation*}
 \xb_{T-2}(\xi_{[T-2]}) + \Xb_{T-1} | \xi_{[T-2]} \geq \yb_{T-2}(\xi_{[T-2]}) + \Yb_{T-1} | \xi_{[T-2]}. 
\end{equation*}
This implies that
\begin{equation*}
 \xb_{T-2}(\xi_{[T-2]}) + \Xb_{T-1} | \xi_{[T-2]} \ssd \yb_{T-2}(\xi_{[T-2]}) + \Yb_{T-1} | \xi_{[T-2]}. 
\end{equation*}
Thus, we infer that \eqref{long-ssd} holds at time $T-2$. Proceeding in the same way, we show the relation for all $t=1,\ldots,T-1$.
\end{proof}

\section{Asset-liability management: problem description and formulation}  \label{sec:ALM}

We consider the decision problem of a financial intermediary managing assets and liabilities exposed to several risk sources over a finite time horizon $T$. The financial intermediary needs to define in which assets to invest at any of the discrete times $t \in {\cal{T}^{\prime}}:=\{0,1,\dots,T-1\}$.
We consider a set $\cal{I}$ of investment assets represented by \textit{exchange traded funds} (ETF), divided into fixed income  ${\cal{I}}_1$, equity ${\cal{I}}_2$ and currencies ${\cal{I}}_3$. The asset universe $\cal{I}^{\prime}$ includes also a cash account labelled as $0$, such as $\cal{I}^{\prime}:= \cal{I} \cup $ $\{ 0 \}$.
At stage $t=0$, we assume that the financial intermediary holds an initial amount $\hat{x}_{i,0}$ of liquid asset $i \in \cal{I}$ and that the proportion of each liquid asset $i \in \cal{I}$ in the portfolio must respect specific lower and upper bounds $\theta^\text{m}_i$ and $\theta^\text{M}_i$ relative to an evolving portfolio value. Fixed income ETFs $i \in {\cal{I}}_1$ are characterized by a constant time to maturity, leading to a deterministic constant duration $\delta_{i,t}^x=\delta_i^x$
for each stage $t$. Assets price returns and gain-loss coefficients are instead considered as random parameters evolving as discrete-time stochastic processes. 

Uncertainty in assets' returns and gain-loss coefficients is represented by means of a non recombining scenario tree, with $\cal{N}$ indicating its set of nodes. For each stage $t \in {\cal T}:= {\cal T}^{\prime} \cup \{T \} $, there is a discrete set of nodes ${\cal N}_t$. The final set ${\cal N}_T$ is the set of nodes called leaves, while the set ${\cal N}_0$ is composed of a unique node, i.e., the root. Each node at stage $t$, except the root, is connected to a unique node at stage $t-1$, which is called ancestor node $a(n)$, and to nodes at stage $t+1$, called successors. For each node $n$ except the leaves (i.e., $n \in {\cal{N}}_t, \ t < T$) there exists a non-empty set of children nodes ${\cal{C}} (n) \in {\cal{N}}_{t+1}$. A scenario is a path through nodes from the root node to a leaf node. We represent with $p_n$ the probability of node $n$ and we have $\sum_{n \in {{\cal{N}}_t}} p_n = 1,\ t \in \cal{T}$.
For each liquid asset $i \in \cal{I}$, we represent by $r_{i,n}$ and $g_{i,n}$ the realization in node $n \in {\cal{N}}_t,t \in \cal{T}$ of price returns and gain-loss coefficients respectively. Details on their evolution will be provided later in Section \ref{sec:assetdyn}.

The other source of uncertainty is given by due payments $L_{j,n}$ generated by liability class $j \in \cal{J}$, reflected in random values in both the liability durations $\delta^\lambda_{j,n}$ and the current fair value of liability contracts $\lambda_{j,n}$. At each stage $t$, the total value of liabilities incurred by the financial intermediary is therefore a stochastic parameter, represented by $\Lambda_n:=\sum_{j \in {\cal J} } \lambda_{j,n}, \ n \in  {\cal{N}}_t$. 
Interest rate immunization is enforced by requiring a narrow mismatching $\overline{\Delta}^{(x,\lambda)}$ between asset and liability durations. 
The core business of the financial intermediary generates over the time horizon uncertain revenues $c_n, \ n \in {\cal N}_t$. Thus, the difference $c_n-\sum_{j \in {\cal J}} L_{j,n}, \ n \in {\cal{N}}_t$ represents the core technical profit in stage $t$.

We assume that in each intermediate stage $t, 1 \leq t \leq T-1$, the financial intermediary, rather than selling assets, may prefer to use a cash overdraft, whose financial cost is determined by an interest rate $r^-_n$ known upon taking the short position.
The continuous variables $b_{n}$ are introduced for this purpose to represent a liquidity credit in node $n \in \mathcal{N}_t$, while the two continuous decision variables $b^+_{n}$ and $b^-_{n}$ correspond to increase or reduction that will concur to determine the intermediary liquidity condition.

In such a framework, the financial intermediary wants to determine a minimal initial investment $K_0$ and an optimal portfolio strategy, by defining the evolution over the time horizon of the assets portfolio and hedge a proportion $\varphi$ of the liabilities, to control the funding ratio. 
Specifically, to represent the amount of each asset $i \in \cal{I}^{\prime}$ held at stage $t \in \cal{T}$, the continuous variables $x_{i,n}, n \in  {\cal{N}}_t$, are introduced. The assets portfolio composition may change over time due to buying and selling decisions that can occur in any but the final stage: $x^+_{i,n}$ and $x^-_{i,n}, n \in  {\cal{N}}_t$, represent buying and selling decisions in node $n$, stage $t \in {\cal T}^{\prime}$. 
Transaction costs associated with buying and selling decisions are represented by parameters $\phi^+$ and $\phi^-$ respectively.


We introduce the following notation:

\textit{Sets}:
\begin{itemize}
\item ${\cal T} = \{ t \ : \ t=0,\dots,T \}$: set of stages;
\item ${\cal T}^{\prime} = \{ t \ : \ t=0,\dots,T-1 \}$: set of stages (last stage excluded);
\item ${\cal N} = \{ n \ : \ n=0,\dots,N \}$: set of nodes of the scenario tree;	
\item ${\cal N}_t \subset {\cal N}$: set of the scenario tree nodes at stage $t \in {\cal T}$;
\item $a(n)$: ancestor of node $n \in {\cal N}_t, \ t \in {\cal T} \setminus \{ 0 \}$;
\item ${\cal C}(n)$: set of children of node $n \in {\cal N}_t, \ t \in {\cal T}^{\prime}$;
\item ${\cal I}^{\prime} = \{ i \ : \ i=0,\dots,I \}$: set of assets (cash included $i=0$);
\item ${\cal I} = \{ i \ : \ i=1,\dots,I \} = \mathcal{I}_1 \cup \mathcal{I}_2 \cup \mathcal{I}_3$: set of investment assets, ETFs;
 \item ${\cal I}_1 = \{ i \ : \ i=1,\dots,I_1 \}\subset \cal{I}$: subset of fixed income assets;
 \item {${\cal I}_2 = \{ i \ : \ i=1,\dots,I_2 \}\subset \cal{I}$: subset of equity assets};
 \item {${\cal I}_3 = \{ i \ : \ i=1,\dots,I_3 \}\subset \cal{I}$: subset of currencies};
\item ${\cal J} = \{ j \ : \ j=1,\dots,J \}$: set of liabilities.

\end{itemize}

\textit{Deterministic Parameters}:
\begin{itemize}
\item $\hat{x}_{i,0}$: initial amount of  asset $i\in\mathcal{I}$ held;
	\item $\theta^{\text{m}}_i $  minimum proportion of $i\in \mathcal{I}$ in the portfolio;
	\item $\theta^{\text{M}}_i$: 	maximum proportion of $i\in \mathcal{I}$ in the portfolio;	
 \item $\delta^x_{i}$:  duration of fixed income asset $i\in \mathcal{I}_1$ in stage $t \in {\cal T}$;
 \item $\overline{\Delta}^{(x,\lambda)}$:  maximum duration mismatching;
 \item 	$\phi^+$:  investment unit transaction cost coefficient;
	\item $\phi^-$:  selling unit transaction cost coefficients; 
\item $\varphi$: funding ratio, or asset-liability ratio, to enforce a funding condition to be satisfied;
\item $q$: maximum proportion of equities in the portfolio; 
 \item $\alpha \in [0,1]$: trade-off parameter  between cumulative debt and investment profits;
 \item $\beta$: penalty coefficient of the initial cash deposit in the objective function.
\end{itemize}

\textit{Stochastic Parameters}:
\begin{itemize}
\item $r_{i,n}$:  price return of asset $i\in\mathcal{I}^{\prime}$ in node $n \in {\cal N}_t, \ t \in {\cal T} \setminus \{ 0 \}$;
\item $g_{i,n}$: gain-loss coefficient of asset $i\in\mathcal{I}$ in node $n \in {\cal N}_t, \ t \in {\cal T} \setminus \{ 0 \}$;
\item $r^-_{n}$: interest rate on debt exposure to the intermediary in node $n \in {\cal N}_t, \ t \in {\cal T} \setminus \{ 0 \}$;
\item $L_{j,n}$: cash outflows associated with liability $j \in {\cal J}$ in node $n \in {\cal N}_t, \ t \in {\cal T} \setminus \{ 0 \}$;
\item $\delta^{\lambda}_{j,n}$:  duration of liability $j\in \mathcal{J}$ in node $n \in {\cal N}$;
\item $\lambda_{j,n}$:  value of liability $j\in \mathcal{J}$ in node $n \in {\cal N}$;
\item $c_n$: core business revenues in node $n \in {\cal N}_t, \ t \in {\cal T} \setminus \{ 0 \}$.
\end{itemize}

\textit{Decision Variables}:
\begin{itemize}
\item $x^+_{i,n}\in\mathbb{R}^{+}$:  amount of asset $i\in \mathcal{I}$ purchased in node $n \in {\cal N}$;
\item $x^-_{i,n}\in\mathbb{R}^{+}$:  amount of asset $i\in \mathcal{I}$ sold in node $n \in {\cal N}$;
\item $x_{i,n}\in\mathbb{R}^{+}$:  amount of asset $i\in \mathcal{I}^{\prime}$ held in node $n \in {\cal N}$;
\item $K_{0} \in \mathbb{R}^+$:  initial investment;
\item $b_{n}\in\mathbb{R}^{+}$: cumulative borrowed money in node $n \in {\cal N}_t, \ t \in {\cal T} \setminus \{ 0 \}$;
\item $b_{n}^+\in\mathbb{R}^{+}$: increase of debt in node $n \in {\cal N}_t, \ t \in {\cal T} \setminus \{ 0 \}$;
\item $b_{n}^-\in\mathbb{R}^{+}$: decrease of debt due to reimbursement in node $n \in {\cal N}_t, \ t \in {\cal T} \setminus \{ 0 \}$.
\end{itemize}
\color{black}

The corresponding optimization model is formulated as follows:
\begin{subequations} 
	\label{eq:NodeALM}
	\begin{align} 
		 \min \,& \rho_0 \circ \dots \circ \rho_{T-1} \left[  \sum_{t \in {\cal T} \setminus \{ 0 \}} \sum_{n \in {\cal N}_t} \left[ 
    \alpha \textcolor{black}{b_{n}} - (1-\alpha) \sum_{i \in {\cal I}} g_{i,n} x^-_{i,n} \right]  \right] + \beta K_0
    \label{eq:obj}\\ \nonumber\\
		\text{ s.t.  } \ & x_{i,0}=\hat{x}_{i,0}+x^+_{i,0}-x^-_{i,0}, \ \  \ \ \ i \in {\cal I},  \label{eq:RootNodeAllocation}\\ 
		& x_{0,0} =
        K_0
        +\sum_{i \in {\cal I}} x^{-}_{i,0} (1-\phi^-) - \sum_{i \in {\cal I}} x^+_{i,0} (1+\phi^+), \label{eq:RootNodeCash} \\
		& x_{i,n}=x_{i,a(n)}(1+r_{i,n})+x^+_{i,n}-x^-_{i,n}, \ \  \ \ \  i \in {\cal I}, \ \ n \in {\cal N}_t, \ \ t \in {\cal T} \setminus \{ 0 \} ,\label{eq:PortfolioRebalancing} \\  
		& x_{0,n}=x_{0,a(n)}(1+r_{0,a(n)}) - \textcolor{black}{b_{a(n)} r_{a(n)}^-}  +  \sum_{i \in {\cal I}} x^-_{i,n}(1-\phi^-) + & \nonumber\\
		&-\sum_{i \in {\cal I}} x^+_{i,n}(1+\phi^+) + c_n-\sum_{j \in {\cal J}} L_{j,n} \textcolor{black}{- b_{n}^- + b_{n}^+}, \ \ n \in {\cal N}_t, \ \ t \in {\cal T} \setminus \{ 0 \}, \label{eq:CashFlows}   \\
 & \textcolor{black}{b_{n} = b_{a(n)}  - b_{n}^- + b_{n}^+, \ \ n \in {\cal N}_t, \ \ t \in {\cal T} \setminus \{ 0 \},}\label{eq:Borrow} \\  
		& \! -\! \sum_{j \in {\cal J}} \! \lambda_{j,n}\overline{\Delta}^{(x,\lambda)} \! \le \! \sum_{i \in {\cal I}_1} \! x_{i,n} \delta^x_{i}\! -\sum_{j \in {\cal J}}\! \lambda_{j,n} \delta^{\lambda}_{j,n} \!\! \le \! \! \sum_{j \in {\cal J}} \! \lambda_{j,n} \overline{\Delta}^{(x,\lambda)}, \ n \in {\cal N}_t, t \in {\cal T}^{\prime},  \label{eq:DurationMatching}  \\
		& \textcolor{black}{ \sum_{i \in {\cal I}^\prime} x_{i,n}(1+r_{i,{\cal C}(n)})  \succeq_{(2)} \varphi \sum_{j \in {\cal J}} \lambda_{j,{\cal C}(n)}, \ \ n \in {\cal N}_{T-1},}   \label{eq:SDorderk} \\
		&\theta^{\text{m}}_i \sum_{i \in {\cal I}} x_{i,n} \le x_{i,n} \le \theta^{\text{M}}_i \sum_{i \in {\cal I}} x_{i,n}, \ \ i \in {\cal I},\ \  n \in {\cal N}_t, \ \ t \in {\cal T}^{\prime},  \label{eq:Diversification}   \\
& \sum_{i \in {\cal I}_2} x_{i,n} \le q \sum_{i \in {\cal I}} x_{i,n}, \ \  n \in {\cal N}_t, \ \ t \in {\cal T}^{\prime},
\label{eq:max_Equities}   \\
 & x_{i,n}^+=x_{i,n}^-=b_{n}^+=
  0, \ \ \ n \in {\cal N}_T.     \label{eq:noterminalrebalance}
	\end{align}
\end{subequations}

Denoting with $\rho_t$ a one-period conditional risk measure, the objective function \eqref{eq:obj} minimizes the initial investment $K_{0}$ with a penalty $\beta$ and, for $t \in {\cal T} \setminus \{ 0 \}$, the compounded cost over the planning horizon of a convex combination via $\alpha \in [0,1]$ between the cumulative debt $b_n$ and the investment profits with negative sign.
The coefficients $\beta$ and $\alpha$ are non-negative fixed inputs, the latter in particular adopted to span alternative risk profiles, from a decision maker focusing only on the liability side (i.e., $\alpha=1$), to one seeking the maximum realized investment profits (i.e., $\alpha=0$).
Notice that the minimization of \( b_n \) in the objective function prevents from excessive borrowing, avoiding arbitrage opportunities when the interest rate on debt exposure is lower than the expected return of some assets in a given node.
Constraint \eqref{eq:RootNodeAllocation} defines the root node optimal investment allocation at  time $t=0$, given an initial input portfolio $\hat{x}_{i,0}$, $i \in {\cal I}$, while equation \eqref{eq:PortfolioRebalancing} represents the rebalancing constraints for asset $i \in {\cal I}$ in subsequent stages. 
Equation \eqref{eq:RootNodeCash} defines the root node cash balance equation in $t=0$ and equation \eqref{eq:CashFlows} traces the evolution of cash surpluses in each stage taking into account interest accrual from previous stage cash balance, payment of interest on cumulative debt $b_{a(n)} r_{a(n)}^-$, selling $x^-_{i,n}$ and buying $x^+_{i,n}$ decisions on the asset portfolio, revenues $c_n$ and costs $\sum_{j \in {\cal J}} L_{j,n}$ associated with the intermediary core business, additional borrowing $b_{n}^+$, and reimboursement $b_{n}^-$. 
Constraint \eqref{eq:Borrow} computes the  debt $b_n$ in node $n$.
Constraint \eqref{eq:DurationMatching} models the assets-liabilities duration matching, which allows a small duration mismatch $\overline{\Delta}^{(x,\lambda)}$ between fixed income assets and liabilities. 
According to Theorem \ref{theorem1}, constraint \eqref{eq:SDorderk} enforces at stage $T-1$ the stochastic dominance of order 2 of the asset portfolio values over the percentage $\varphi$ of the liability portfolio on the children nodes. A $\varphi$ equal to 1 would reflect a \textit{fully-funded} liability, whereas we would have under-funding for $\varphi<1$ and over-funding for $\varphi>1$. 
As established by Lemma 5.27 and Lemma 5.32 in \cite{DDARriskbook}, equation \eqref{eq:SDorderk} defines a convex constraint.
Constraint \eqref{eq:Diversification}, as customary in institutional ALM, helps enforcing portfolio diversification by imposing minimum and maximum investment proportions.
Finally, constraint \eqref{eq:max_Equities} imposes an upper bound to the proportion of equities in the portfolio, and \eqref{eq:noterminalrebalance} rules out possible assets purchases or liquidation and debt increase at the end of the planning horizon.

\section{Uncertainty model} \label{sec:statmodel}

The ALM model implementation requires the specification of a rich set of random coefficients, assumed, in this setting, to follow a discrete non recombining tree process. We present in the following subsections the stochastic models adopted to derive the full set of coefficients of problem \eqref{eq:NodeALM}. The following economic and financial risk sources are accounted for and motivate the model specification based on investment horizon $T$ and a liability evaluation horizon $T_\lambda$: 
\begin{itemize}
    \item The fluctuations of the term structure of interest rates and inflation have a joint impact on the asset portfolio and the liability of the intermediary through their pricing and duration mismatching.
    \item Credit risk fluctuations in the economy have an impact on corporate returns and exogenous liability costs.
    \item Every asset class has specific risk factors driving their future behaviour.
    \item All the above jointly determine the risk exposure of the intermediary.
\end{itemize}

\subsection{Yield curve, inflation and credit spread models} \label{sec:economics}

Asset returns and liability costs depend on the evolution of inflation and the term structure of interest rates. For the latter we have implemented the popular \textit{Nelson-Siegel-Svensson} model in the dynamic arbitrage-free version proposed by \cite{Chris2009}. Let $y_{t,\tau}$ denote the yield quoted at time $t$ over the term $\tau$. We have:
\begin{equation} \label{eq:NSmodel}
y_{t,\tau}=b^y_{1,t}+b^y_{2,t} e^{-\frac{\tau}{\gamma_t}}+b^y_{3,t} \frac{\tau}{\gamma_t} e^{-\frac{\tau}{\gamma_t}}.
\end{equation}
Equation \eqref{eq:NSmodel} relies on a 3-factor model with factors reflecting level, slope and convexity of the curve, as functions of $b^y_{1,t}, \ b^y_{2,t}$ and $b^y_{3,t}$, from which a long-term yield $b^y_{1,t}$ and an instantaneous yield $b^y_{1,t}+b^y_{2,t}$ for the instantaneous short rate dynamics $y_{t,0}=y_t$ can be derived. The parameter $\gamma_t$ represents a decay factor, here expressed and estimated as a linear function of the coefficients $b_{j,t}^y$, $j=1,2,3$:
\begin{equation} \label{eq:LambdaModel}
\gamma_t=a^y_0 + a^y_{1} b^y_{1,t} + a^y_{2}b^y_{2,t} + a^y_{3}b^y_{3,t} + \epsilon_t^y,
\end{equation}
where $a^y_j, j=0,1,2,3$, are the coefficient processes and $\epsilon_t^y$ are the residuals, here supposed to be normally distributed and correlated. We have adopted multivariate Ordinary Least Squares (OLS) estimation and calibrated the model to enforce arbitrage free conditions, following \cite{Chris2009}.

We consider as reference for the inflation process $\pi_t$ the annual Consumer Price Index (CPI) dynamics, specified as a simple mean-reverting model with the long-term mean set at the European Central Bank (ECB) $2\%$ target. Given an initial state $\pi_0$, for $t \in {\cal T}$ with monthly increments $\Delta t$, we assume a dynamic:

\begin{equation} \label{eq:CPI}
\pi_t=\pi_{t-\Delta t}+a^{\pi}(0.02-\pi_{t-\Delta t})\Delta t+\sigma^{\pi}\sqrt{\pi_{t-\Delta t}}\sqrt{\Delta t}  \epsilon^{\pi}_t,
\end{equation}
where $\epsilon_t^\pi$ are the residuals that we assume normally distributed and correlated, while parameters $a^{\pi}$ and $\sigma^{\pi}$ are the coefficient processes to be estimated on the data history.

Jointly with the yield curve, the spread model for Investment Grade (IG) borrowers in the euro area defines an explanatory variable for corporate costs as specified next. Credit spread $s_t^{IG}$ random dynamics are described by the following autoregressive model: 
\begin{equation} \label{eq:spread}
s_t^{IG}=c_0+c_1 s^{IG}_{t-\Delta t}+c_2 y_{t}+\epsilon^{s}_t,
\end{equation}
where $c_j, j=0,1,2$, are the coefficient processes to be estimated on data, $\epsilon_t^s$ are the residuals, supposed to be normally distributed and correlated and $y_t$ is the instantaneous short rate. The fitted term structure of interest rates and the credit spread process lead jointly to the definition of the borrowing costs faced by the intermediary, with $r^-_t=y_{t,1}+s_t^{IG}$ both expressed per annum and defined over the 1 year term.

We refer to Appendix \ref{Stat} for the estimation of the coefficients associated with the statistical models introduced in this paragraph.

\subsection{Liability model} \label{sec:liabdyn}

We consider a liability process of a representative insurance intermediary with an \textit{investment grade} credit rating that over time funds its activity relying on incoming premiums generated by claims and life contracts. 
The modeling framework is inspired by \cite{consigli2011b}, where the intermediary is classified as an investment-grade BBB borrower in terms of credit rating based on current market conditions.
Following the approach presented in \cite{consigli2011a}, we introduce a simplified liability model, where revenues and costs evolve according to a simple linear stochastic process. Specifically, given initial long-term estimates $c_0$ and $L_{j,0}$, cash inflows and revenues $c_t$ at time $t$ and insurance costs and compensations to underwriters $L_{j,t}$ for liability class $j$ at time $t$ are described by the following equations:

\begin{eqnarray} \label{eq:netrevenues}
c_{t+1}&=c_t(1+\rho_{t+1}),  \\
\label{eq:rev_inc} \rho_t&=\mu_{\rho}+\sigma_{\rho} \epsilon^{\rho}_{t},  \\
\label{eq:netpay} L_{j,t+1}&=L_{j,t}(1+\xi_{j,t+1}),   \\
\label{eq:pay_inc} \xi_{j,t}&=\mu_{j,\xi}+\sigma_{j,\xi} \epsilon^{\xi}_{j,t}, \end{eqnarray}
where $\mu_{\rho}$ ($\mu_{j,\xi}$) is the average increase for cash inflows (outflows), $\sigma_{\rho}$ ($\sigma_{j,\xi}$) is the volatility of cash inflows (outflows), and $\epsilon^{\rho}_{t}$ ($\epsilon^{\xi}_{j,t}$) are the normally distributed residuals.
By subjectively modifying the premiums and claims distributions, we can generate \textit{stressed} liability scenarios and verify the impact on the optimal funding conditions and investment policy.
Our primary interest in this application is twofold: (i) to assess how the introduction of stressed liability scenarios impacts the optimal investment and funding strategy, and (ii) to further evaluate the effectiveness of multistage stochastic dominance constraints within the ALM framework. 

Given the liability cash flows $L_{j,t}$, the current liability obligation $\lambda_{j,t}$ is defined as the discounted value of the expected payments as follows: 
\begin{equation}  \label{eq:lambda1}
\lambda_{j,t}=\mathbb{E}_{t}\left[\sum_{h=t}^{t + T_{\lambda} } e^{-y_{t,h}(h-t)} L_{j,h}  \right].
\end{equation}
We estimate the nodal liability values as scenario dependent discounted cash flows in descending nodes: at each node this estimation includes all cash flows projected over $T_{\lambda}$ years.
The associated duration $\delta^\lambda_{j,t}$ is defined as follows: 
\begin{equation}  
\label{eq:duration_liab}
\delta^\lambda_{j,t}=\mathbb{E}_{t}\left[\frac{\sum_{h=t}^{t+ T_{\lambda}} (h-t) e^{-y_{t,h}(h-t)} L_{j,h} } {\lambda_{j,t}} \right]  .
\end{equation}
Discounting is then attained by backward recursion along the tree relying on the nodal realizations of the yield curve from the model described above.

This simplified approach enables us to assess the impact of different liability scenarios within the ALM model while avoiding unnecessary model complexity.
For more refined liability valuation, alternative methods such as the replicating portfolio approach, the least squares Monte Carlo method, or actuarial models (see \cite{black1973pricing}, \cite{pelsser2016difference}, \cite{dickson2016insurance}, and references therein) could be employed, depending on the level of detail required.

\subsection{Asset and currency returns} \label{sec:assetdyn}
{We consider a partition of the asset universe ${\cal I}$ into three classes: for fixed income ETFs ${\cal I}_1:=\{i=1,\dots, I_1\}$, for equity ETFs ${\cal I}_2:=\{i=I_1 +1,\ldots,I_2\}$ and for currencies ${\cal I}_3:=\{i=I_2 +1,\ldots,I_3\}$. The set ${\cal I}_1$ is further partitioned into treasury fixed income assets $i=1,\dots,I_4,\  I_4 < I_1,$ and corporate fixed income assets $i=I_4+1,\dots,I_1$.}

From observed data, asset price returns are computed as $r_{i,t}=\frac{v_{i,t}}{v_{i,t-\Delta t}}-1$, 
with $v_{i,t}$ denoting the ETF value at time $t$. 
In order to investigate the relationship between asset returns and selected explanatory variables, we use linear regression models tailored to each asset class.
For treasury fixed income assets,
we have:
\begin{equation} \label{eq:A1model}
    r_{i,t}=b_{i,0}+b_{i,1} r_{i,t-1}+b_{i,2} y_{t-1,\delta_{i}^x}+b_{i,3} \pi_t+\epsilon_{i,t}, \quad i=1,\dots,I_4 ,
\end{equation}
where $y_{t-1,\delta_{i}^x}$ is the yield to maturity at time $t-1$ of asset $i$ with underlying duration $\delta_{i}^x$, $\pi_t$ is the EU inflation rate at time $t$ given in \eqref{eq:CPI} and $\epsilon_{i,t}$ are the residuals. Ordinary least squares estimation is employed to determine the regression coefficients $b_i=(b_{i,0}, b_{i,1}, b_{i,2}, b_{i,3})^\top, \ i=1,\dots,I_4$.

{For corporate assets we consider the following linear regression model:
\begin{equation} \label{eq:A2model}
r_{i,t}=b_{i,0}+b_{i,1} r_{i,t-1}+b_{i,2} y_{t-1,\delta_{i}^x}+b_{i,3} s^{IG}_{t}+b_{i,4} r_{I_1+1,t} + \epsilon_{i,t}, \quad i=I_4+1,\dots,I_1,
\end{equation}
which  assumes a dependence on the credit spread variable $s_t^{IG}$ and on small cap returns $r_{I_1+1,t}$ of equity $I_1+1 \in \mathcal{I}_2$.  The vector of regression coefficients to be estimated is given by $b_i=(b_{i,0}, b_{i,1}, b_{i,2}, b_{i,3}, b_{i,4})^\top$.}

{Equity asset returns for  $i \in {\cal I}_2$ are considered for large and small caps in the equity market plus emerging markets. These assets become of primary importance in the long-term when trying to stochastically dominate the liability costs that may increase significantly upon increasing uncertainty. The correspondent autoregressive model is:
\begin{equation} \label{eq:A3model}
r_{i,t}=b_{i,0}+b_{i,1} r_{i,t-1}+b_{i,2} y_{t,1} + b_{i,3} \pi_t + b_{i,4} f_{t}+ \epsilon_{i,t}, \quad i=I_1+1,\dots,I_2.
\end{equation}
Model (\ref{eq:A3model}) postulates dependence of equity returns on previous returns, on the 1-year interest rate $y_{t,1}$, on the inflation rate $\pi_t$ and on the \textit{term spread} $f_t$ between the 10- and the 1-year interest rates, whose relevance is related to its ability to capture long-term economic expectations.} The corresponding vector of regression coefficients is $b_i=(b_{i,0}, b_{i,1}, b_{i,2}, b_{i,3}, b_{i,4})^\top, \ i=I_1+1,\ldots,I_2$. 

The autoregressive model for currency returns, for $i \in \mathcal{I}_3$, is:
\begin{equation} \label{eq:A4model}
    r_{i,t}=b_{i,0}+b_{i,1} r_{i,t-1}+b_{i,2} y_{t,0.25}+b_{i,3} \pi_t+\epsilon_{i,t}, \quad i=I_2 + 1,\dots,I_3.
\end{equation}
Model (\ref{eq:A4model}) postulates dependence of currency returns on previous returns, on the 3-month interest rate $y_{t,0.25}$ and on the inflation rate $\pi_t$. The corresponding vector of regression coefficients is $b_i=(b_{i,0}, b_{i,1}, b_{i,2}, b_{i,3})^\top, \ i=I_2+1,\ldots,I_3$. 
We refer to Appendix \ref{Stat} for the results of the estimation of the regression coefficients $b_i=(b_{i,0}, b_{i,1}, b_{i,2}, b_{i,3}, b_{i,4})^\top, \ i\in\mathcal{I}$.

Once defined the asset returns, we can compute the gain and loss coefficients $g_{i,t}$, which are associated with selling decisions of any type of asset. For $h \in {\cal T}, i \in {\cal I}$, let $\rho_{i,h}:=\Pi_{s=1}^{h}(1+r_{i,s})-1$. We define the average gain and loss coefficient per unit selling at time $t$ as:

\begin{equation} \label{eq:gaincoeff}
g_{i,t}:=\frac{1}{t} \sum_{h=1}^t \rho_{i,h}.
\end{equation}
We refer to  Appendix \ref{App_Scen} for the scenario generation algorithm adopted to derive the tree processes of the entire set of stochastic parameters.

\section{Solution method}  \label{sec:decomposition}

The idea is to adopt the dynamic programming formulation of the multi-stage problem and to solve it recursively by considering a new version of the multi-cut method, see \cite{Rus2003Stoc}, in which additional event cuts approximate the stochastic order constraints and further cuts as in \cite{gulten2015two} approximate the risk measure in the objective function. The objective is in form \eqref{eq:obj}, which represents a time-consistent dynamic risk measure. 
 

The problem can be solved recursively as follows.
At the last stage, we calculate for every leaf node $n\in \mathcal{N}_T$ the cumulative debt:  
\begin{align}
Q_{n,T}=\min\; & \alpha b_n \notag\\
\text{s.t. } 
  &  x_{i,n}=x_{i,a(n)}(1+r_{i,n}), \notag\\
         & x_{0,n}=x_{0,a(n)}(1+r_{0,a(n)}) -b_{a(n)} r^-_{a(n)} +c_n -\sum_{j \in \mathcal{J}} L_{j,n}, \notag  \\
        & b_n=b_{a(n)}  - b^-_n. \notag       
\end{align}
For the nodes $n\in\mathcal{N}_{T-1}$, we calculate:
\begin{align}
Q_{n,T-1} = \min\; & \alpha b_n - (1-\alpha) \sum_{i \in \mathcal{I}} g_{i,n} x^-_{i,n}  +  \varrho_n[\mathcal{Q}_{T}| n]\notag\\
\text{s.t. }\; 
  &  x_{i,n}=x_{i,a(n)}(1+r_{i,n})+x^+_{i,n}-x^-_{i,n}, \ \  i \in \mathcal{I}, \label{eq:n_PortfolioRebalancing} \\
  & x_{0,n}=x_{0,a(n)}(1+r_{0,a(n)}) -b_{a(n)}r^-_{a(n)} +\sum_{i \in {\cal I}} x^-_{i,n}(1-\phi^-) + & \nonumber\\
 &-\sum_{i \in \mathcal{I}} x^+_{i,n}(1+\phi^+) + {c_n-\sum_{j \in {\cal J}} L_{j,n}} + b^+_n - b_n^-, \;\label{eq:n_CashFlows}   \\ 
 & b_n=b_{a(n)} + b^+_n - b^-_n, \\
 & \theta^{\text{m}}_i \sum_{i \in \mathcal{I}} x_{i,n} \le x_{i,n} \le \theta^{\text{M}}_i \sum_{i \in \mathcal{I}} x_{i,n}, \ \ i \in \mathcal{I},\;  
 \label{eq:n_Diversification} \\
 & {\sum_{i \in \mathcal{I}_2} x_{i,n} \le q \sum_{i \in \mathcal{I}} x_{i,n},   
 \label{eq:n_tot_equities} }\\ 
& -\sum_{j \in \mathcal{J}} \lambda_{j,n} \overline{\Delta}^{(x,\lambda)} \le \sum_{i \in \mathcal{I}_1} x_{i,n} \delta^x_{i}-\sum_{j \in \mathcal{J}} \lambda_{j,n} \delta^{\lambda}_{j,n} \le \sum_{j \in \mathcal{J}} \lambda_{j,n} \overline{\Delta}^{(x,\lambda)}, \label{eq:n_DurationMatching} \\
  & \sum_{i \in {\cal I}^\prime} x_{i,n}(1+r_{i,{\cal C}(n)})  \succeq_{(2)} \varphi \sum_{j \in {\cal J}} \lambda_{j,{\cal C}(n)}. \;  \label{eq:n_SDorderk}
\end{align}
For the nodes $n\in\mathcal{N}_{t}$, $t=1\dots,T-2$, we calculate:
\begin{align}
Q_{n,t} = \min\; & \alpha b_n - (1-\alpha) \sum_{i \in \mathcal{I}} g_{i,n} x^-_{i,n}  +  \varrho_n[\mathcal{Q}_{t+1}| n]\notag\\
\text{s.t. }\; 
  & \eqref{eq:n_PortfolioRebalancing}-\eqref{eq:n_DurationMatching}. \notag 
\end{align}
For the root node:
\begin{align*}
\min\;\; & \beta 
K_0 +\varrho_0[\mathcal{Q}_1]\\
\text{s.t. }\; 
 & x_{i,0}=\hat{x}_{i,0}+x^+_{i,0}-x^-_{i,0}, \ \  i \in {\cal I}, \\ 
        & x_{0,0}=
        K_0+\sum_{i \in {\cal I}} x^{-}_{i,0} (1-\phi^-) - \sum_{i \in {\cal I}} x^+_{i,0} (1+\phi^+),\\
 & \eqref{eq:n_Diversification} - \eqref{eq:n_DurationMatching}.
\end{align*}



We start by solving the problems at the leaf nodes with some initial guess for the variables at their ancestor nodes.
Assume that we carry out iteration $\ell$. 
For each problem at node $n\in\mathcal{N}_t$, we obtain its optimal value, denoted 
$\bar{v}_n^\ell$, the optimal solution, denoted $({x}_n^\ell,({x}_n^+)^\ell,({x}_n^-)^\ell,{b}_n^\ell,({b}_n^+)^\ell,({b}_n^-)^\ell)$, and the optimal Lagrange multipliers $d_n^\ell\in\mathbb{R}^{|\mathcal{I}|+2}$ associated with constraints about re-balancing of the cash, assets, and the cumulative debt. This information provides an objective cut at the ancestor node of $n$ of form:
\begin{align*}
v_{a(n)}\;& \geq \bar{v}_n^\ell + \langle - T_n^\top d_n^\ell, (x_{a(n)},b_{a(n)}) - ({x}_{a(n)}^\ell,{b}_{a(n)}^\ell)\rangle \\
& = -\langle T_n^\top d_n^\ell, (x_{a(n)},b_{a(n)})\rangle + \alpha_n^\ell,\quad \text{with }
\alpha_n^\ell = \bar{v}_n^\ell + \langle T_n^\top d_n^\ell, ({x}_{a(n)}^\ell,{b}_{a(n)}^\ell)\rangle.
\end{align*}
Here $T_n$ is the matrix containing the coefficients associated with the ancestor variables of node $n$, with elements $d_{00} = 1+ r_{0,a(n)}$, $d_{ii} = 1+ r_{i,n}$ for $i\in\mathcal{I}$, $d_{ii}= 1$ for $i= |\mathcal{I}|+2$, and $d_{0i}= r_{a(n)}^{-}$ for $i= |\mathcal{I}|+2$.
We shall gather the objective cuts for the objective function of node $n$ constructed until iteration $\ell$ in the set $J_o^\ell(n)$.

Furthermore, at node $n$, having solved the problems for all successor nodes $m\in\mathcal{C}(n)$, we solve an auxiliary problem 
\begin{equation}
\label{p:risk_approximation}
\max_{\mu\in\mathcal{A}_\varrho} \sum_{m\in\mathcal{C}(n)} p_{n,m}\mu_m \bar{v}_m^\ell,
\end{equation}
where $\mathcal{A}_\varrho$ is the convex subdifferential  $\varrho[0]$ in the dual representation of the risk measure $\varrho.$ Let $\mu_n^\ell$ be the solution of that problem. As it is a subgradient of the risk measure, it provides a cut in the approximating problem for node $n$ of the following form:
\[
w_n\geq \langle\mu_n^\ell,v\rangle.
\]
These cuts approximating the risk function at node $n$ that are constructed until iteration $\ell$ are gathered in the set $J_r^\ell(n)$. We also need a parameter $\underline{w}$ to impose a lower bound on the value of the risk measure $\varrho_n.$

The ordering constraint is approximated according to the quantile method presented in \cite{DenMar2012}, see also \cite{dentcheva2010inverse}. This means that at node  $n\in\mathcal{N}_t$, we compare the random variable $\Lambda_n$ with realizations $\Lambda_{n,m}=\sum_{j\in\mathcal{J}} \lambda_{j,m}$, $m\in\mathcal{C}(n)$, and the random variable $X_n$ with realizations 
$X_{n,m}=\sum_{i\in\mathcal{I}^\prime} x_{i,n}(1+r_{i,m})$, where $r_{i,m}$ is associated with node $m\in\mathcal{C}(n).$

In order to impose the stochastic dominance constraints \eqref{eq:SDorderk} assuming the order $k=2$, we use the following method.
We denote $S^1= \{1,\dots, |\mathcal{C}(n)|\}.$
\begin{description}
\item {\bf Algorithm to impose stochastic dominance} 
\item[Step 0:] Set $\iota =1$, $J_e^{\iota}(n)=\{ S^1\}$, and $X_{n,m}^1=\sum_{i\in\mathcal{I}^\prime} x_{i,n}^\ell(1+r_{i,m})$ for all $m\in\mathcal{C}(n)$.
\item[Step 1:] Solve the problem:
\begin{equation}
\label{cp1i}
\begin{aligned}
\min \;&\; \alpha b_n - (1-\alpha) \sum_{i \in \mathcal{I}} g_{i,n}x_{i,n}^{-} +w_n + \beta 
K_0\\
\text{s.t. } &\; w_n\geq \langle\mu_n^j,v\rangle \quad j\in J_r^\ell(n), \; v \in \mathbb{R}^{|\mathcal{C}(n)|},\\
&\; v_m\geq -\langle T_m^\top d_m^j, (x_n,b_n)\rangle + \alpha_m^j,\quad j\in J_o^\ell (m),\; m\in\mathcal{C}(n),\\
&\; x_{i,n}=x_{i,{a(n)}}^{\ell-1}(1+r_{i,n})+x^+_{i,n}-x^-_{i,n}, \ \  i \in \mathcal{I}, \\
  &\; x_{0,n}=x_{0,a(n)}^{\ell-1}(1+r_{0,a(n)})-b_{a(n)}^{\ell-1}r^-_{a(n)}+\sum_{i \in \mathcal{I}} x^-_{i,n}(1-\phi^-) + \\
 &\quad-\sum_{i \in \mathcal{I}} x^+_{i,n}(1+\phi^+) + c_n -\sum_{j \in \mathcal{J}} L_{j,n} + b^+_n - b^-_n, \\ 
 &\;b_n= b_{a(n)}^{\ell -1} + b^+_{n} - b^-_{n} ,  \\
 &\theta^{\text{m}}_i \sum_{i \in \mathcal{I}} x_{i,n} \le x_{i,n} \le \theta^{\text{M}}_i \sum_{i \in \mathcal{I}} x_{i,n}, \ \ i \in \mathcal{I},\\ 
& { \sum_{i \in \mathcal{I}_2} x_{i,n} \le q \sum_{i \in \mathcal{I}} x_{i,n},  
} \\
 &\; - \sum_{j \in \mathcal{J}} \lambda_{j,n} \overline{\Delta}^{(x,\lambda)} \le \sum_{i \in \mathcal{I}_1} x_{i,n} \delta^x_{i}-\sum_{j \in \mathcal{J}} \lambda_{j,n} \delta^{\lambda}_{j,n} \le \sum_{j \in \mathcal{J}} \lambda_{j,n} \overline{\Delta}^{(x,\lambda)},\\
&\; \frac{1}{P(S^j)} \sum_{m\in S^j} p_{n,m} X_{n,m}^\iota \geq\frac{1}{P(S^j)} F^{(-2)}(\Lambda_n;P(S^j)),\quad S^j\in J_e^\iota(n).
 \end{aligned}
\end{equation}
Let $X_n^{\iota}$ be the new random variable associated with the solution of problem \eqref{cp1i}.
If $n \notin \mathcal{N}_{T-1}$, then index the solutions of problem \eqref{cp1i} by $\ell$ and stop. Otherwise, continue.
\item[Step 2:] Consider the sets $A^{\iota}_\eta = \{ X_n^{\iota} \le \eta\}$ and let
\begin{equation}
\label{deltai}
\delta_{\iota} = \sup_\eta \Big\{  \frac{1}{P(A^{\iota}_\eta)} F^{(-2)}(\Lambda_n;P(A^{\iota}_\eta)) - \frac{1}{P(A^{\iota}_\eta)}\sum_{m\in A^{\iota}_\eta} p_{n,m} X_{n,m}^\iota  : P(A^{\iota}_\eta) >0\Big\}.
\end{equation}
If $\delta_{\iota} \le 0$, then index the solutions of problem \eqref{cp1i} by $\ell$ and stop. Otherwise, continue.
\item[Step 3:]   Find $\eta^{\iota}$ such that $P(X_n^{\iota} \le \eta^{\iota}) >0$ as well as
\begin{equation}
\label{newcut1}
\frac{1}{P(A^{\iota}_{\eta^\iota})} \sum_{m\in A^{\iota}_{\eta^\iota}} p_{n,m} X_{n,m}^\iota  - \frac{1}{P(A^{\iota}_{\eta^\iota})} F^{(-2)}(\Lambda_n;P(A^{\iota}_{\eta^{\iota}})) \le - \frac{\delta_{\iota}}{2},
\end{equation}
are satisfied.
\item[Step 4:]
Set $S^{\iota} = A^{\iota}_{\eta^{\iota}}$, $J_e^{\iota +1}(n)= J_e^\iota(n)\cup \{ S^{\iota}\}$, increase $\iota$ by one, and go to Step 1.
\end{description}

The solution of the problem $v_n^\ell$ provides a lower bound for the recourse function
$Q_n(x_{a(n)},b_{a(n)})$, while $w_n^\ell$ is a lower bound for the risk measure associated with node $n.$

If the problem is infeasible, we can construct a feasibility cut,
\begin{equation}
\label{3.5m}
\gamma_n^\ell + \langle \tilde{d}_n^\ell, (x_{a(n)},b_{a(n)})\rangle \leq 0.
\end{equation}
The feasibility cuts remain valid for the true cost-to-go function.

We refer to the approximate problem \eqref{cp1i} at each node of the scenario tree as
$\mathcal{P}(n)$.
Each of the problems $\mathcal{P}(n)$ maintains and updates the following data:
its current solution $({x}_n,({x}_n^+),({x}_n^-),{b}_n$, $({b}_n^+),({b}_n^-))$, convex polyhedral models of the cost-to-go
functions $\underline{Q}^{(j)}(\cdot)$ of
its successors $m\in \mathcal{C}(n)$  (if any), and the current approximation
$v_n$ and $w_n$ of the optimal value of its own cost-to-go function and the risk measure at $n$.
The operation of each subproblem is as follows:

\begin{description}
\item[Step 1.]
If $n$ is not the root node, retrieve from the ancestor problem $\mathcal{P}(a(n))$ its current approximate solution $(x_{a(n)},b_{a(n)})$.
\item[Step 2.]
If $n$ is not a leaf node,  retrieve from each successor problem $\mathcal{P}(m)$,  $m\in \mathcal{C}(n)$,
all new objective and feasibility cuts and update the approximations
of their cost-to-go functions $\underline{Q}_n(\cdot)$. Update the approximation of its risk measure by solving problem \eqref{p:risk_approximation}.
\item[Step 3.]
Solve the problem \eqref{cp1i}.
\item[(a)]
If it is solvable, update its solution and its optimal value. If $n$ is not the root node and $v_n$
increased, construct a new objective cut.
\item[(b)]
If the problem is infeasible, and $n$ is not the root node,
construct a new feasibility cut. If $n$ is the root node,
then stop, because the entire problem is infeasible.
\item[Step 4.]
Wait for the command to activate again, and  then go to Step 1.
\end{description}

It remains to describe the way in which these subproblems are initiated,
activated in the course of the solution procedure,
and terminated.
We assume that we know a sufficiently large number $M$ such that
each cost-to-go function can be bounded from below by $-M$. Our
initial approximations of the successors' functions are just
\[
\underline{Q}^{(j)}(\cdot) = -M.
\]
At the beginning, no ancestor solutions are available, but we can initiate
each subproblem with some arbitrary point $(x_{a(n)},b_{a(n)})$.

There is much freedom in determining the order
in which the subproblems are solved. Three rules have to be observed.
\begin{enumerate}
\item
There is no sense to activate a subproblem $\mathcal{P}(n)$ whose ancestor's solution
did not change, and whose successor problems $\mathcal{P}(m)$, $m\in\mathcal{C}(n)$,
 did not generate any new cuts since this problem was activated last.
\item
If a subproblem $\mathcal{P}(n)$ has a new solution, each of its successors
$\mathcal{P}(m)$, $m\in\mathcal{C}(n)$ has to be activated some time after this solution has been obtained.
\item
If a subproblem $\mathcal{P}(n)$ generates a new cut, i.e., if
it is infeasible or has a new optimal value $v_n$,
its ancestor $\mathcal{P}(a(n))$ has to be activated some time after this cut has been generated.
\end{enumerate}
We shall terminate the method if Rule 1 applies to all subproblems, in which case
we claim that the current solutions constitute the optimal solution
of the entire problem. The other stopping test is the infeasibility test at Step 3(a)
for the root node. It is obvious, because we operate with relaxations here, and if
a relaxation is infeasible, so is the true problem.

Now we argue that the method discovers infeasibility of the problem or converges to a solution of it. 
If the method stops because the Rule 1 applies to all subproblems, i.e., no subproblem needs to be activated, then, we claim that the current solutions constitute the optimal policy of the entire problem. If the method stops because of infeasibility at the root node, then the whole problem is infeasible,  because we operate with relaxations. Hence, if a relaxation is infeasible, so is the true problem. 
We also observe that for each set of decisions, the algorithm employed to impose the SD constraint terminates in finitely many steps discovering infeasibility or identifying an optimal solution. This is due to the fact that the Lorenz functions of random variables with finitely many realizations are piece-wise linear; cf. also \cite{DenMar2012}, Theorem 4. The approximation of the risk measure for a fixed random variable converges to its true value due to the convergence of the cutting plane method because the subdifferential set  $\mathcal{A}_\varrho$ is a closed and bounded convex set. In the case of the mean-semideviation of first order, or the average value at risk combined with the expected value, we shall obtain an exact calculation after finitely many steps because the dual set is polyhedral, otherwise, we need to terminate the approximation when a prescribed numerical accuracy is reached. In such a case, only finitely many cuts are used to approximate the measure of risk up to the prescribed accuracy. Finally, the multicut method approximates the optimal value of the recourse function by objective cuts and its domain by feasibility cuts. It is convergent due to the convergence of the cutting plane method. For polyhedral functions, the method converges in finitely many iterations (see \cite{Rus2003Stoc}, Chapter 3).   

\section{Computational evidence} \label{sec:computations}

In this section, an extended set of computational results is presented with the aim of validating the proposed methodology and discussing the most relevant financial evidence. 
Following the ALM problem in \eqref{eq:NodeALM}, we consider an ALM manager seeking a minimal initial capital injection, sufficient however to fund an investment strategy, with periodic revision, able to cover all liabilities over the following 5 years.
The risk associated with total debt and investment losses is controlled by means of the mean-semideviation risk measure.
We present through the section the results collected assuming either a \textit{base liability scenario}, we may also refer to as \textit{ongoing ALM} scenario, or a \textit{stressed liability scenario}, as the one recently experienced in insurance markets. 
In particular, the main features of the data set used to generate the assets' and liability scenario trees are first summarized with their statistical properties in Section \ref{sec:datainput}. We then present in Section \ref{sec:method_validation} the  evidence collected on the decomposition method developed to solve the optimization problem. Section \ref{sec:rootpflio} and Section \ref{sec:intraterisk} focus respectively on the results on risk preferences and interest rates exposures induced by the optimal solutions. The impact of stochastic dominance constraints is analyzed specifically in Section \ref{sec:SDimpact} with final results on the intermediary solvency conditions over a 5 year planning horizon. 

\subsection{Data inputs and experimental design} \label{sec:datainput}

We take the perspective of a generic European insurance intermediary with a 5-year planning horizon for strategic asset allocation and liability hedging, see \cite{consigli2012b}.
The asset universe includes the following \textit{ Exchange Traded Funds} (ETF), or benchmarks (in round brackets the ID ticker for \textit{Yahoo! Finance}, see \cite{Yahoo}): 
\begin{itemize}
    \item Money market index: \textit{UCITS ETF C-EUR} (SMART.MI);
    \item 1-3 year bond index: \textit{iShares Govt Bond 1-3yr UCITS ETF} (IBGS.L);
    \item 5-7 year bond index: \textit{Xtrackers II Eurozone Govt Bond 5-7 UCITS ETF} (DBXR.DE);
    \item 10 year bond index: SPDR Bloomberg 10+ Year Euro Govt Bond UCITS ETF (SYBV.DE);
    \item IG corporate bond index: \textit{iShares iBoxx Investment Grade Corporate Bond ETF} (LQD);
    \item Intermediate-term bond index: \textit{Vanguard Total Bond Market Index Fund} (BND);
    \item 3 year inflation linked bond index: \textit{FlexShares iBoxx 3-Year Target Duration TIPS Index Fund} (TDTT);
    \item 5 year inflation linked bond index: \textit{FlexShares iBoxx 5-Year Target Duration TIPS Index Fund} (TDTF);
    \item Large cap equity index: \textit{iShares Core MSCI Europe UCITS ETF EUR} (IMEU.AS);
     \item Small cap equity index: \textit{iShares Russell 2000 ETF} (IWM);
     \item Emerging markets equity index: \textit{iShares MSCI Emerging Markets ETF} (EEM);
     \item Japan market equity index: \textit{iShares MSCI Japan ETF} (EWJ);
     \item US market equity index: \textit{iShares MSCI USA SRI UCITS ETF USD} (SUAU.AS);
     \item EUR-USD exchange rate: \textit{EUR/USD} (EURUSD).
\end{itemize}

The statistical models are calibrated with a data history of monthly observations from December 2018 to December 2022. The subset  ${\cal I}_1$ consists of eight fixed income bonds for duration matching (SMART.MI, IBGS.L, DBXR.DE, SYBV.DE, LQD, BND, TDTT, TDTF).  The subset ${\cal I}_2$ includes a global equity ETF (IMEU.AS), the small cap equity ETF (IWM), an ETF for emerging markets (EEM), and equity ETFs for Japan (EWJ) and US (SUAU.AS). The subset ${\cal I}_3$ includes the EUR-USD exchange rate.
The decision space is thus ${\cal I}={\cal I}_1 \cup {\cal I}_2 \cup {\cal I}_3$. 
Following the details in Table \ref{duration}, the ETFs in class ${\cal I}_1$ are classified as \textit{constant-to-maturity} (CTM) fixed income benchmarks carrying by construction a relatively stable duration coefficient. 
The asset-liability duration matching will rely on such simplifying assumptions. Insurance liabilities do instead carry a scenario dependent duration, as explained below.

\begin{table}[h!]
\begin{small} 
\begin{tabular}{|lll|} \hline
Index $i \in {\cal I}_1$ & ETF code & Duration $\delta_{i}^x$ \\ \hline
Money market index & Smart.MI  & $0.25$ year \\
1-3 year bond index & IBGS.L & $1.6$ years \\
5-7 year bond index & DBXR.DE & $5.6$ years\\
10 year bond index  & SYBV.DE  & $14.5$ years  \\
IG corporate bond index & LQD  & $8.4$ years \\
Intermediate-term bond index & BND & 6.1 years \\
3 year inflation linked bond index & TDTT & 3 years \\
5 year inflation linked bond index & TDTF & 5 years \\ 
\hline  
\end{tabular}
\caption{Durations $\delta_{i}^x$ of fixed-income assets $i \in {\cal I}_1$.}
\label{duration}
\end{small}
\end{table}

The following settings are assumed in the case study: 
\begin{itemize}
\item ALM planning horizon of $5$ years with stages ${\cal T}:=\{0, 1, 2, 3, 5\}$ years.
\item A scenario tree with branching $10^4$ resulting into ${\cal N}_4=10000$ scenarios and $N=11111$ nodes. Stochastic values for asset returns and liability costs are determined by applying the models described in Section \ref{sec:statmodel} and the scenario generation procedure summarized in Appendix \ref{App_Scen}. 
\item Different risk-reward trade-offs in the objective function can be considered by varying the coefficients $\alpha$ and $\beta$ in equation \eqref{eq:obj}. Parameter $\alpha$ determines the trade-off between the cumulative debt and the investment profit. By contrast, $\beta$ acts uniquely on the initial portfolio estimate and will determine the optimal initial risk capital to be allocated in the form of a given root node portfolio allocation. 
In our numerical experiments, we assume $\alpha=\{0.25,0.5,0.75\}$ and $\beta=\{1,2\}$. 
\item Liabilities estimated by projecting cash outflows over $T_{\lambda}=5$ years.
\item Initial levels $\hat{x}_{i,0}=0$ and investment lower and upper bounds $\theta_i^{\text{m}}=0$ and $\theta_i^{\text{M}}=0.3$ respectively on all assets $i \in {\cal I}$. 
Additionally, we analyse the impact of asset diversification by considering a case where a minimum investment constraint of $\theta_i^{\text{m}}=3\%$ on all assets is imposed.
Transaction costs are set to $\phi^{+}=\phi^{-}=0.001$. 
\item A maximum percentage of equity assets in the portfolio $q=0.4$. We also analyse the case $q=1$ to assess the impact of this constraint on the model’s outputs.
\item Maximum duration mismatching between assets and liabilities $\overline{\Delta}^{(x,\lambda)}=\{6,3\}$ months. These values help analysing the effectiveness of the portfolio immunization strategies. The ALM manager will allocate the investment portfolio to hedge against liability interest rate exposure.
\item Weight for the funding ratio $\varphi=\{0,0.8,1,1.1\}$, used to analyse the impact of varying liability funding requirements on portfolio strategy.
\end{itemize}

\paragraph*{Asset statistics and scenario tree process}

Table \ref{r_Scen} compares the average monthly returns and standard deviation from historical data to those associated at the horizon with the tree process for each asset. 

\begin{table}
    \centering
    \begin{tabular}{l|l|rrrrrr}
    \hline
 Asset & Series & Mean& Std & Skewness& Kurtosis& 5\%-Q & 1\%-Q
\\
\hline
  \multirow{2}{*}{Smart.MI} &History& --0.0003& 0.0002& 1.0244& 4.0505& --0.0005& -- 0.0005
\\
  &Simulation& --0.0002& 0.0005& 0.0076& 3.1370& --0.0082& --0.0116
\\
\hline
  \multirow{2}{*}{IBGS.L}  &History&  0.0022&  0.0207&  0.4366&  4.3537&  --0.0345&  --0.0542
\\
          &Simulation&  0.0032&  0.0238&  0.3925&  3.9203&  --0.0021&  --0.0478
\\
\hline
 \multirow{2}{*}{DBXR.DE} &History&  0.0013&  0.0073&  --0.2766&  2.5383&  --0.0129&  --0.0153
\\
          &Simulation&  0.0016&  0.0089&  --0.0868&  2.9967&  --0.0201&  --0.0270
\\
\hline
 \multirow{2}{*}{SYBV.DE} &History&  0.0012&  0.0235&  --0.2725&  2.8639&  --0.0450&  --0.0590
\\
          &Simulation&  0.0013&  0.0172&  0.0188&  2.9555&  --0.0309&  --0.0428
\\
\hline
  \multirow{2}{*}{LQD} &History&  0.0019&  0.0186&  --0.8626&  8.0608&  --0.0246&  --0.0838
\\
          &Simulation&  0.0021&  0.0149&  --0.0124&  3.1017&  --0.0302&  --0.0816
\\
\hline
  \multirow{2}{*}{BND} &History&  --0.0002&  0.0121&  --0.1971&  5.0905&  --0.0206&  --0.0421
\\
          &Simulation&  --0.0006&  0.0249&  --0.4010&  3.2838&  --0.0180&  --0.0638
\\
\hline
 \multirow{2}{*}{TDTT} &History& 0.0009& 0.0062& --0.0498& 2.8686& --0.0080& --0.0184
\\
  &Simulation& 0.0006& 0.0064& --0.0250& 2.9646& --0.0143& --0.0187
\\
\hline
 \multirow{2}{*}{TDTF} &History& 0.0007& 0.0071& --0.0135& 3.2622& --0.0120& --0.0284
\\
  &Simulation& 0.0007& 0.0067& --0.0138& 3.3240& --0.0122& --0.0266
\\
\hline
 \multirow{2}{*}{EEM} &History& 0.0038& 0.0493& -0.6429& 4.9842& --0.0669& --0.1824
\\
  &Simulation& 0.0040& 0.0536& -0.3416& 4.3776& --0.0641& --0.2002
\\
\hline
 \multirow{2}{*}{IMEU.AS} &History& 0.0037& 0.0449& --0.4644& 5.8507& --0.0621& --0.1698
\\
  &Simulation& 0.0040& 0.0562& --0.4293& 3.6215& --0.0787& --0.1597
\\
\hline
  \multirow{2}{*}{IWM} &History& 0.0096& 0.0637& --0.5466& 9.3377& --0.0885& --0.2850
\\
  &Simulation& 0.0091& 0.0536& 0.3595& 4.6976& --0.0763& --0.2829
\\
\hline
  \multirow{2}{*}{EWJ} &History& 0.0022& 0.0577& --0.2953& 3.6475& --0.0828& --0.1200
\\
  &Simulation& 0.0022& 0.0572& --0.1170& 2.9680& --0.0887& --0.1090
\\
\hline
  \multirow{2}{*}{SUAU.AS} &History & 0.0017 & 0.0579 & 0.2261& 4.2296& --0.0601& 
--0.1438\\
  &Simulation & 0.0016 &  0.0568 & 0.2010& 3.1189& --0.0555& 
--0.1511\\
\hline
  \multirow{2}{*}{EURUSD} &History& --0.0003& 0.0268& --0.1637& 4.7531& --0.0473& --0.0760
\\
  &Simulation& --0.0004& 0.0266& --0.1193& 3.1136& --0.0209& --0.0634
\\
\hline
    \end{tabular}
    \caption{\label{r_Scen} Mean, standard deviation (Std), skewness, kurtosis, 5\%-quantile (5\%-Q), 1\%-quantile (1\%-Q) of monthly returns $r_{i,t}$ of assets $i \in {\cal I}$: historical versus simulation evidence at the end of year 5.}
    \label{tab:stats_ass}
\end{table}

\paragraph*{Liability scenarios}

\begin{figure} [ht!]
    \centering
    \includegraphics[width=0.9\linewidth]{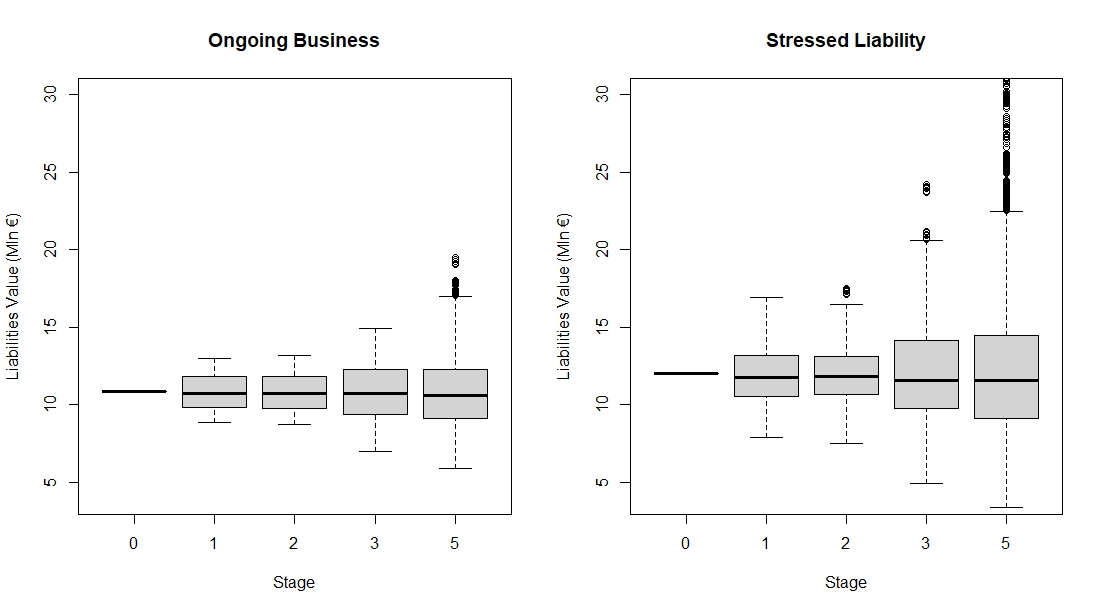}   
    \caption{Liability estimates $\Lambda_n=\sum_{j \in {\cal J}} \lambda_{j,n}, n \in {\cal N}_t, t\in{\cal T}$ in the \textit{ongoing business} scenario (left) and in the \textit{stressed liability} scenario (right).}
    \label{fig:Distr_Liab}
\end{figure}

The liability estimates follow the case study developed in \cite{consigli2011b} for a large P\&C company with an estimated first year cash inflows due to collected premiums of 4.2 million \texteuro{},  mean $\mu_{\rho}=0.5\%$, volatility $\sigma_{\rho}=1\%$, and first year cash outflows due to casualties associated with a liability $L_{1,n}$ of 2.2 million \texteuro{}, an estimated $1\%$ average annual increase $\mu_{1,\xi}$ and a $3\%$ volatility $\sigma_{1,\xi}$. These forecasts correspond to the base scenario case.
We consider a stressed scenario by assuming an annual average increase of $\mu_{1,\xi}=5\%$ and a $5\%$ volatility $\sigma_{1,\xi}$. 
Based on these estimates we derive two possible evolutions of the insurance liability reserves $\Lambda_n = \sum_{j \in {\cal J}} \lambda_{j,n}$, shown in Figure \ref{fig:Distr_Liab}: in particular under the base scenario the current (time 0) liability estimate is $\Lambda_0=10.21$ million \texteuro{}, while under the stressed scenario, \textit{ceteris paribus}, this amount increases to $\Lambda_0=11.82$ million \texteuro{}.

The introduction of a stressed liability scenario serves two purposes: first, to evaluate the financial intermediary’s solvency under more severe liability conditions; second, to further validate the adopted SSD methodology and problem formulation.

\subsection{Problem decomposition: numerical results} 
\label{sec:method_validation}


On the available data set we applied the solution method described in Section \ref{sec:decomposition}. We chose the mean-semideviation of order 1 as a risk measure $\rho_t$ (see \cite{Rus2003Stoc}), which allowed as to formulate problems \eqref{cp1i} as linear programs. Specifically, we set the weighting parameter between mean and semideviations in the risk-measure specification to $0.1$ 
As a lower bound on the value of the risk measure we considered $\underline{w}=-10^6$. Parameters $\alpha$ and $\beta$ in the objective function \eqref{eq:obj} have been set to $0.5$ and $1$ respectively. Solutions for different values of these parameters are investigated in the following paragraph. All computational experiments were run on an ASUS laptop with a 3 GHz Intel Core i7-5500U Processor and 4 GB of RAM using solver Gurobi under GAMS 24.7.4 environment. 
We start by solving the problems at the leaf nodes by considering the solution of the worst-case liability scenario as initial guess for the
variables at their ancestor nodes $(x^0_{a(n)},b^0_{a(n)})$. 
This choice allows us to start the numerical procedure with a feasible solution in the leaf nodes.
One iteration of the algorithm consists in solving all $N=11111$ problems ${\cal P}(n)$ in the scenario tree.  
Indeed, as explained in Section \ref{sec:decomposition}, the proposed solution algorithm decomposes the problem into a sequence of two-stage problems, where each two-stage problem is associated with a specific node in the scenario tree. 
As a result, the size of each subproblem can vary depending on the stage of the corresponding node. 
Table \ref{Sol Time} provides further details on the size of each subproblem.

At each iteration, problems ${\cal P}(n), \ n \in {\cal N}$, are solved following a backward approach, from the leaf nodes to the root. 
In terms of order in which the subproblems are solved, to take into account the three rules described in Section \ref{sec:decomposition}, we control the set of problems to be solved by means of a binary parameter $ON(n)$, which is set to $1$ if problem ${\cal P}(n)$ needs to be solved and to $0$ otherwise.
Once problem ${\cal P}(n)$ is solved, we set $ON(n)=0$ and we compare the solution provided at the current iteration with the previous one: if problem ${\cal P}(n)$ has a new solution, we send an activation signal to the ancestor and to the children nodes of node $n$ (i.e., $ON_{a(n)}=1$ and $ON_m=1, \ m \in {{\cal C}(n)}$). At each iteration, we only solve problems with $ON_{n}=1$, since $ON_{n}=0$ implies that the solution of the problem ${\cal P}(n)$ does not differ from the one determined before.
The number of iterations required by the algorithm described above in the base case instance is 7. Thanks to the specific initialization procedure, problems ${\cal P}(n), \ n \in {\cal N}$, are always feasible in all iterations: this prevents from the generation of feasibility cuts. The total CPU time needed to solve the problem is 4147 seconds, corresponding to $1$ hour $9$ minutes and $7$ seconds. 
With respect to the stressed liability case, this stressed situation does not affect the performance of the solution algorithm, keeping the number of iterations the same as in the base case, not requiring further feasibility cuts.
We present in Table \ref{Sol Time} additional evidences on the solution times 
at the last iteration of the solution algorithm. 
\begin{table}
\begin{small}
\begin{tabular} {|lccccc| } \hline
Stage $t \in {\cal T}$ & $0$ & $1$ & $2$ & $3$ & $4$ \\ \hline
N. of problems ${\cal P}(n), \ n \in {\cal N}_t$ & 1 & 10 & 100 & 1000 & 10000 \\ 
N. of ${\cal P}(n)$ continuos variables & 44 & 43 & 43 & 43 & 35 \\
N. of ${\cal P}(n)$ constraints & 50 & 46 & 46 & 56 & 14 \\
CPU time spent on stage $t$& $0.02\%$ & $1.01\%$ & $2.05\%$ & $9.82\%$ & $87.10\%$ \\ \hline
\end{tabular}
\caption{\label{Sol Time}Number and size of problems ${\cal P}_n, n \in {\cal N}_t$ and CPU time allocation over stages $t \in {\cal T}$.}
\end{small}
\end{table}
As can be noticed, despite the smaller size of problems in the leaf nodes, due to their large cardinality, most of the CPU time is spent solving the problem at the last stage.
Results on stochastic dominance constraints will be described later in Section \ref{sec:SDimpact}.

\subsection{Risk preferences, optimal debt exposure and initial capital allocation} \label{sec:rootpflio}

The coefficient $\alpha \in [0,1]$ defines a trade-off between cumulative debt and investment profit. The coefficient $\beta$ represents instead a penalty coefficient on the initial investment.
We can then associate different risk profiles of the ALM manager to each pair $\{\alpha,\beta\}$. 
A $\{0.5,1\}$-type of ALM manager would assign same relevance to the debt minimization and investment profits maximization, while being sufficiently safe with the initial capital invested. 
A $\{0.75,1\}$-type of manager would instead be more concerned with debt limitation, while a $\{0.5,2\}$-type  would focus on the minimization of the initial investment. 
In this paragraph we discuss the impact on model results of different weights assigned to parameters $\alpha$ and $\beta$.
Specifically, Tables \ref{Debt} and \ref{Prof} show respectively the average across the scenarios and the standard deviation of the debt position and of cumulative investment profits in each stage $t \in {\cal T}\setminus\{0\}$. 
Regarding the investment profits, as enforced by equation \eqref{eq:noterminalrebalance}, the financial intermediary is restricted from making any asset purchases or liquidations at the terminal stage. As a result, cumulative investment profits remain constant from year 3 to year 5. For this reason, the cumulative profits for year 5 are not presented in Table \ref{Prof}, as they are identical to those for year 3.
Table \ref{Root Comp} provides the optimal initial investment value $K_0$, the funding ratio at stage 0, which is the asset to liability ratio at the beginning of the planning period, and the average funding ratio at the horizon.

\begin{table}[ht!]
\begin{small}
\begin{tabular}{ |m{0.2cm} | m{0.4cm} 
                m{0.4cm} 
                m{0.4cm} 
                m{0.4cm} 
                m{0.4cm}
                m{0.4cm}  
                m{0.4cm}  
                m{0.8cm} | 
                m{0.8cm} 
                m{0.8cm} 
                m{0.8cm} 
                m{0.8cm} | }\hline
& & & & & & & & & \multicolumn{4}{c|}{Year}\\ 
\# & $\mu_{1,\xi}$ & $\sigma_{1,\xi}$ & $\alpha$ & $\beta$ & $q$ & $\theta_i^\text{m}$ & $\varphi$ & $\overline{\Delta}^{(x,\lambda)}$  & $1$ & $2$ & $3$ & $5$\\ \hline

1 & $1\%$ & $3\%$ & $0.5$ & $1$ & 0.4 & 0 & 1 & 0.5 & 5.98 (3.22) & 6.40 (4.82) & 0.13 (0.59) & 0.00 (0.00)
\\ \hline \hline

2 & $1\%$ & $3\%$ & $0.25$ & $1$ & 0.4 & 0 & 1 & 0.5 & 8.41 (3.11) & 9.26 (2.54) & 0.20 (0.80) & 0.00 (0.00)
\\  \hline
3 & $1\%$ & $3\%$ & $0.75$ & $1$ & 0.4 & 0 & 1 & 0.5 & 1.98 (0.52) & 0.70 (2.56) & 0.03 (0.57) & 0.00 (0.00) 
\\ \hline
\hline

4 & $1\%$ & $3\%$ & $0.5$ & $2$ & 0.4 & 0 & 1 & 0.5 & 9.24 (1.01) & 6.51 (4.68) & 0.19 (0.74) & 0.00 (0.00) 
\\ \hline \hline

5 & $5\%$ & $5\%$ & $0.5$ & $1$ & 0.4 & 0 & 1 & 0.5 & 6.30 (5.22)& 7.78 (4.79)& 0.33 (0.77)& 0.00 (0.00)\\ \hline \hline

6 & $1\%$ & $3\%$ & $0.5$ & $1$ & 1 & 0 & 1 & 0.5 & 4.81 (2.28) & 5.03 (3.94) & 0.29 (0.85) & 0.00 (0.00)
\\ \hline \hline

7 & $1\%$ & $3\%$ & $0.5$ & $1$ & 0.4 & 3\% & 1 & 0.5 & 5.01 (1.50) & 2.29 (2.52) & 0.69 (1.10) & 0.00 (0.00)
\\ \hline \hline

8 & $5\%$ & $5\%$ & $0.5$ & $1$ & 0.4 & 0 & 1 & 0.25 & 6.36 (3.21)& 7.80 (4.82)& 0.41 (0.53)& 0.00 (0.00)\\ \hline \hline

9 & $1\%$ & $3\%$ & $0.5$ & $1$ & 0.4 & 0 & 0 & 0.5 & 4.30 (4.42) & 5.31 (4.74) & 0.00 (0.00) & 0.00 (0.00)
\\ \hline 

10 & $5\%$ & $5\%$ & $0.5$ & $1$ & 0.4 & 0 & 0 & 0.5 & 4.45 (4.21)& 6.02 (4.83)& 0.00 (0.00)& 0.00 (0.00)\\ \hline \hline

11 & $1\%$ & $3\%$ & $0.5$ & $1$ & 0.4 & 0 & 0.8 & 0.5 & 5.54 (3.01) & 6.25 (4.81) & 0.08 (0.38) & 0.00 (0.00)
\\ \hline 

12 & $5\%$ & $5\%$ & $0.5$ & $1$ & 0.4 & 0 & 0.8 & 0.5 & 6.12 (4.81)& 7.01 (4.64)& 0.16 (0.83)& 0.00 (0.00)\\ \hline \hline

13 & $1\%$ & $3\%$ & $0.5$ & $1$ & 0.4 & 0 & 1.1 & 0.5 & 6.21 (4.01) & 6.50 (4.85) & 0.33 (1.18) & 0.22 (0.47)
\\ \hline 

14 & $5\%$ & $5\%$ & $0.5$ & $1$ & 0.4 & 0 & 1.1 & 0.5 & 6.45 (4.74)& 7.11 (4.81)& 0.55 (3.35)& 0.41 (0.21)\\ \hline 
\end{tabular}
\caption{\label{Debt} Debt exposure in each stage $t \in {\cal T}\setminus\{0\}$: average amount over the scenarios in million of euros and, in parenthesis, standard deviation for different values of parameters $\mu_{1,\xi}$, $\sigma_{1,\xi}$, $\alpha$, $\beta$, $q$, $\theta_i^\text{m}$, $\varphi$ and $\overline{\Delta}^{(x,\lambda)}$.}
\end{small}
\end{table}

\begin{table}[ht!]
\begin{small}
\begin{tabular}{ |m{0.2cm} | m{0.4cm} 
                m{0.4cm} 
                m{0.4cm} 
                m{0.4cm} 
                m{0.4cm}
                m{0.4cm}  
                m{0.4cm}  
                m{0.8cm} | 
                m{0.8cm} 
                m{0.8cm}  
                m{0.8cm} | }\hline
& & & & & & & & & \multicolumn{3}{c|}{Year}\\ 
\# & $\mu_{1,\xi}$ & $\sigma_{1,\xi}$ & $\alpha$ & $\beta$ & $q$ & $\theta_i^\text{m}$ & $\varphi$ & $\overline{\Delta}^{(x,\lambda)}$  & $1$ & $2$ & $3$ \\ \hline

1 & $1\%$ & $3\%$ & $0.5$ & $1$ & 0.4 & 0 & 1 & 0.5 & $-0.61$ (0.46) & $-1.11$ (1.63) & 5.15 (4.00)
\\ \hline \hline

2 & $1\%$ & $3\%$ & $0.25$ & $1$ & 0.4 & 0 & 1 & 0.5 & $-0.95$ (0.59) & $-0.34$ (0.56) & 7.23 (2.71)
\\  \hline
3 & $1\%$ & $3\%$ & $0.75$ & $1$ & 0.4 & 0 & 1 & 0.5 & $-0.52$ (0.52) & $-0.54$ (1.70) & 4.03 (3.17)
\\ \hline
\hline

4 & $1\%$ & $3\%$ & $0.5$ & $2$ & 0.4 & 0 & 1 & 0.5 & $-0.73$ (0.40) & $-0.62$ (1.07) & 3.57 (2.85)
\\ \hline \hline

5 & $5\%$ & $5\%$ & $0.5$ & $1$ & 0.4 & 0 & 1 & 0.5 & $-0.79$ (0.43) & $-0.96$ (1.14) & 3.95 (3.86)
\\ \hline \hline 

6 & $1\%$ & $3\%$ & $0.5$ & $1$ & 1 & 0 & 1 & 0.5 & $-1.34$ (0.81) & $-1.19$ (1.04) & 6.32 (5.35) 
\\ \hline \hline

7 & $1\%$ & $3\%$ & $0.5$ & $1$ & 0.4 & 3\% & 1 & 0.5 & $-0.33$ (0.31) & $-0.29$ (1.03) & 1.59 (1.53) 
\\ \hline \hline

8 & $5\%$ & $5\%$ & $0.5$ & $1$ & 0.4 & 0 & 1 & 0.25 & $-0.86$ (0.40) & $-1.01$ (0.74) & 3.81 (3.97)
\\ \hline \hline

9 & $1\%$ & $3\%$ & $0.5$ & $1$ & 0.4 & 0 & 0 & 0.5 & $-0.99$ (0.59) & $-0.10$ (0.50) & 7.44 (2.83)
\\ \hline

10 & $5\%$ & $5\%$ & $0.5$ & $1$ & 0.4 & 0 & 0 & 0.5 & $-0.88$ (0.78) & $-0.07$ (0.51) & 4.59 (2.53)
\\ \hline \hline

11 & $1\%$ & $3\%$ & $0.5$ & $1$ & 0.4 & 0 & 0.8 & 0.5 & $-0.86$ (0.48) & $-0.42$ (0.70) & 6.96 (2.83) 
\\ \hline 

12 & $5\%$ & $5\%$ & $0.5$ & $1$ & 0.4 & 0 & 0.8 & 0.5 & $-0.64$ (0.38) & $-0.84$ (1.01) & 4.31 (2.57) 
\\ \hline \hline

13 & $1\%$ & $3\%$ & $0.5$ & $1$ & 0.4 & 0 & 1.1 & 0.5 & $-0.61$ (0.46) & $-1.06$ (1.58) & 5.04 (4.06)
\\ \hline 

14 & $5\%$ & $5\%$ & $0.5$ & $1$ & 0.4 & 0 & 1.1 & 0.5 & $-0.98$ (0.69) & $-1.01$ (0.94) & 0.70 (5.69)
\\ \hline 
\end{tabular}
\caption{\label{Prof} Cumulative investment profits in each stage $t \in {\cal T}^{\prime}\setminus\{0\}$: average amount over the scenarios in million of euros and, in parenthesis, standard deviation for different values of parameters $\mu_{1,\xi}$, $\sigma_{1,\xi}$, $\alpha$, $\beta$, $q$, $\theta_i^\text{m}$, $\varphi$ and $\overline{\Delta}^{(x,\lambda)}$.}
\end{small}
\end{table}

Lines 1, 2 and 3 of Tables \ref{Debt}, \ref{Prof} and \ref{Root Comp} are associated with the solution of the ongoing business scenario with different relevance assigned to the debt minimization and investment profits maximization. 
As can be noted, the three instances show a similar evolution of the cumulative debt, which increases in the early stages and then decreases at year 3. At the horizon, the debt is fully reimbursed in all three instances.
In terms of profits, the three instances show a similar behaviour, with investment profit being mainly negative in the first stage and subject to a sharp increase in year 3 due to the liquidation of part of the asset portfolio.
From Tables \ref{Debt} and \ref{Prof} it can be noticed how higher values of $\alpha$ decrease the average profits and reduce the average debt, especially in the early stages: from a financial perspective this is consistent with the extra-weight assigned to the debt minimization in the objective function of the problem.
%
\begin{table}[ht!]
\begin{small}
\begin{tabular}{ |m{0.2cm} | m{0.4cm} 
                m{0.4cm} 
                m{0.4cm} 
                m{0.4cm} 
                m{0.4cm}
                m{0.4cm}  
                m{0.4cm}  
                m{0.8cm} | 
                m{0.8cm} 
                m{0.8cm}  
                m{0.8cm} | }\hline
\# & $\mu_{1,\xi}$ & $\sigma_{1,\xi}$ & $\alpha$ & $\beta$ & $q$ & $\theta_i^\text{m}$ & $\varphi$ & $\overline{\Delta}^{(x,\lambda)}$  & $K_0$ &  $FR_0$ & $FR_5$ \\ \hline

1 & $1\%$ & $3\%$ & $0.5$ & $1$ & 0.4 & 0 & 1 & 0.5 & 15.779 & 1.545 & 1.156 
\\ \hline \hline

2 & $1\%$ & $3\%$ & $0.25$ & $1$ & 0.4 & 0 & 1 & 0.5 & 15.159 & 1.485 & 1.529  
\\  \hline
3 & $1\%$ & $3\%$ & $0.75$ & $1$ & 0.4 & 0 & 1 & 0.5 & 16.785 & 1.642 & 1.196  
\\ \hline
\hline

4 & $1\%$ & $3\%$ & $0.5$ & $2$ & 0.4 & 0 & 1 & 0.5 & 10.106 & 0.990 & 1.016  
\\ \hline \hline

5 & $5\%$ & $5\%$ & $0.5$ & $1$ & 0.4 & 0 & 1 & 0.5 & 18.660 & 1.579 & 1.221   
\\ \hline \hline 

6 & $1\%$ & $3\%$ & $0.5$ & $1$ & 1 & 0 & 1 & 0.5 & 11.610 & 1.137 & 1.959   
\\ \hline

7 & $1\%$ & $3\%$ & $0.5$ & $1$ & 0.4 & 3\% & 1 & 0.5 & 16.170 & 1.584 & 1.386  
\\ \hline \hline

8 & $5\%$ & $5\%$ & $0.5$ & $1$ & 0.4 & 0 & 1 & 0.25 & 19.611 & 1.659 & 1.111
\\ \hline \hline

9 & $1\%$ & $3\%$ & $0.5$ & $1$ & 0.4 & 0 & 0 & 0.5 & 9.095 & 0.891 & 0.797  
\\ \hline

10 & $5\%$ & $5\%$ & $0.5$ & $1$ & 0.4 & 0 & 0 & 0.5 & 7.748 & 0.656 & 0.725 
\\ \hline \hline

11 & $1\%$ & $3\%$ & $0.5$ & $1$ & 0.4 & 0 & 0.8 & 0.5 & 10.106 & 0.990 & 1.007 
\\ \hline 

12 & $5\%$ & $5\%$ & $0.5$ & $1$ & 0.4 & 0 & 0.8 & 0.5 & 9.298 & 0.787 & 1.130  
\\ \hline \hline

13 & $1\%$ & $3\%$ & $0.5$ & $1$ & 0.4 & 0 & 1.1 & 0.5 & 17.752 & 1.739 & 1.310  
\\ \hline 

14 & $5\%$ & $5\%$ & $0.5$ & $1$ & 0.4 & 0 & 1.1 & 0.5 & 20.821 & 1.761 & 1.746 
\\ \hline 

\end{tabular}
\caption{\label{Root Comp} Optimal initial investment value $K_0$ (in million of euros), portfolio allocation, funding ratio at stage 0 ($FR_0$) and average funding ratio at year 5 ($FR_5$) for different values of parameters $\mu_{1,\xi}$, $\sigma_{1,\xi}$, $\alpha$, $\beta$, $q$, $\theta_i^\text{m}$, $\varphi$ and $\overline{\Delta}^{(x,\lambda)}$.}
\end{small}
\end{table}
From Table \ref{Root Comp} it can be noted that the different weights assigned to $\alpha$ have a limited impact on the initial investment.
Additionally, the evolution of the funding ratio indicates that, in all three instances, the optimal portfolio strategy ensures that assets remain sufficient to cover liabilities throughout the planning horizon, thereby preserving long-term financial stability.

With regard to parameter $\beta$, line 4 of Tables \ref{Debt}, \ref{Prof} and \ref{Root Comp} assesses the impact of a higher weight $\beta=2$ for the initial investment in the ongoing business scenario. 
With respect to the case $\beta=1$ (line 1), results show a significant reduced initial investment $K_0$, which leads to lower investment profits and the need to increase the cumulative debt to increase the liquidity to effectively manage the portfolio, especially in the early stages.

\begin{figure}[ht!]
    \centering
    \begin{minipage}{1\textwidth}
        \centering
        \includegraphics[width=\textwidth]{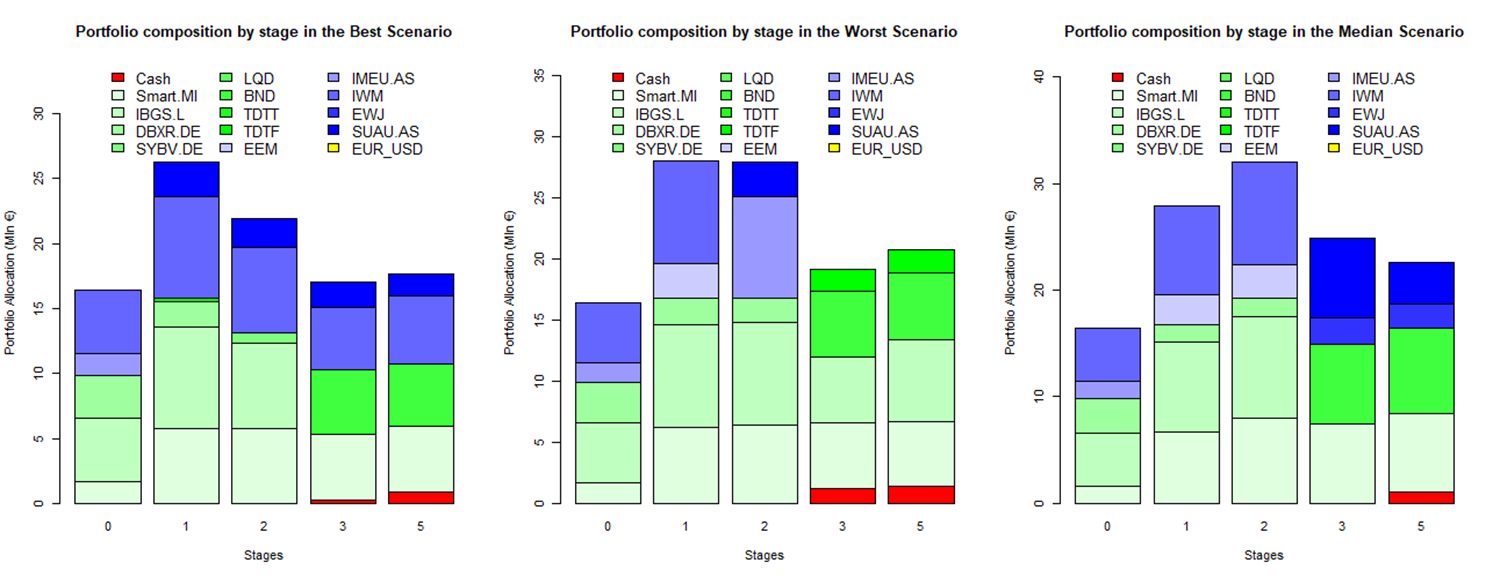}
        \subcaption{Ongoing Business}\label{fig:panel1}
    \end{minipage}
     \vspace{0.5cm}
    \begin{minipage}{1\textwidth}
        \centering
        \includegraphics[width=\textwidth]{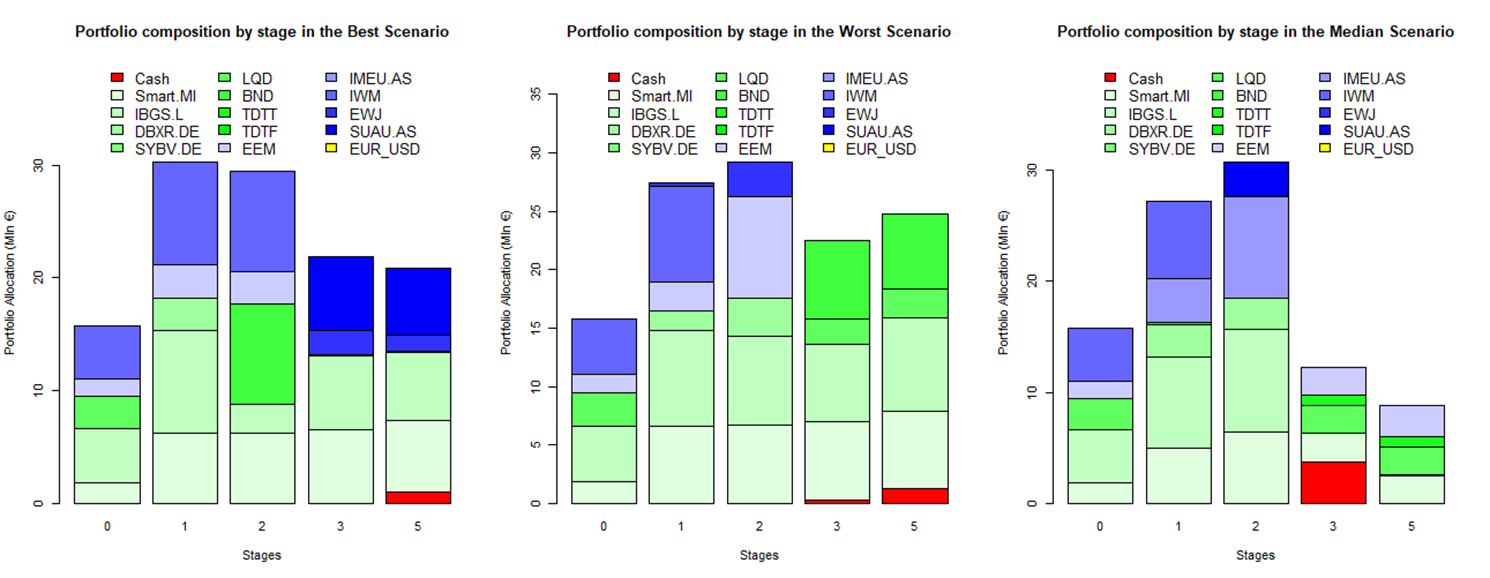}
        \subcaption{Stressed Liability}\label{fig:panel2}
    \end{minipage}
    \caption{Portfolio evolution over time for selected scenarios in both the ongoing business and stressed liability cases. Left to right: best, median and worst scenarios.}
    \label{fig:Portfolio}
\end{figure}

Line 5 is associated with the stressed liability scenario.
Compared to the ongoing business case, the higher values and volatility of liabilities lead to poorer portfolio performances. 
Moreover, as can be noted in Table \ref{Root Comp}, the stressed condition for liabilities results in a significant increase of the capital to be allocated at $t=0$: due to the high cost of capital, such a scenario would be highly undesirable.
To provide an example of the portfolio’s evolution over time, we graphically represent in Figure \ref{fig:Portfolio} the portfolio’s composition in the two liability cases across three selected scenarios: the scenario with minimum liabilities at the horizon (i.e., the \textit{best scenario}), the scenario with maximum liabilities at the horizon (i.e., the \textit{worst scenario}), and the scenario with intermediate liabilities at the horizon (i.e., the \textit{median scenario}).
We observe a highly differentiated portfolio composition across the scenarios. In particular, the terminal allocation between fixed income (represented in green in Figure \ref{fig:Portfolio}) and equities (represented in blue) varies significantly depending on the scenario.
For example, in the worst scenario, starting from year 3, the portfolio is composed exclusively of cash and fixed income assets, which are required to match the duration of liabilities. 
In contrast, both the best and median scenarios display a more balanced portfolio, maintaining a significant equity component throughout the planning horizon.
Despite these differences in asset allocation, we observe that the equity allocation reaches the maximum allowed threshold of 40\% in all stages in best and median scenarios, highlighting that the equity constraint may play a substantial role in shaping the overall portfolio composition.
From Figure \ref{fig:Portfolio}, we can also observe that the portfolio typically consists of a relatively small number of assets, with around 5 to 6 different assets across the scenarios.
These observations on portfolio composition motivate the experiments reported in lines 6 and 7 of Tables \ref{Debt}, \ref{Prof}, and \ref{Root Comp}, where we analyse portfolio performance under different composition requirements. 
Specifically, we investigate the impact of removing the equity constraint (line 6) and enforcing greater asset diversification by imposing for each asset a minimum percentage in the portfolio $\theta_i^\text{m}=3\%$ (line 7).
In the absence of a maximum equity investment constraint, the optimal solution involves an initial investment of 11.610 million \texteuro{}, composed of 64.7\% equities and 35.3\% fixed income. 
Despite the reduction in initial investment, the possibility to manage the portfolio with greater flexibility leads to a significant improvement in portfolio performances, generating higher investment profits while reducing average debt. 
On the other hand, the portfolio diversification requirement results in a substantial deterioration in the financial intermediary's performance, with the optimal portfolio management requiring a higher initial investment, which, however, corresponds to a modest profit.

\subsection{Control of interest rate risk} \label{sec:intraterisk}

The ALM model \eqref{eq:NodeALM} includes a specific set of constraints associated with the exposure to yield curve fluctuations. The ALM manager seeks an optimal strategy while imposing a relatively strict constraint on duration mismatching between assets and liabilities. 
We are then considering only first-order impact of yield curve movements along the tree jointly on the asset and the liability portfolios. 
When asset duration exceeds liability duration, increasing interest rates will negatively affect the portfolio, whereas if liabilities have a longer duration than assets, decreasing interest rates will be detrimental.
To capture these asymmetric effects of interest rate movements, we assess interest rate risk exposure through the duration-matching constraint \eqref{eq:DurationMatching}, 
distinguishing between positive and negative mismatches.

\begin{figure}[ht!]
    \centering    \includegraphics[width=1\linewidth,height=0.4\linewidth]{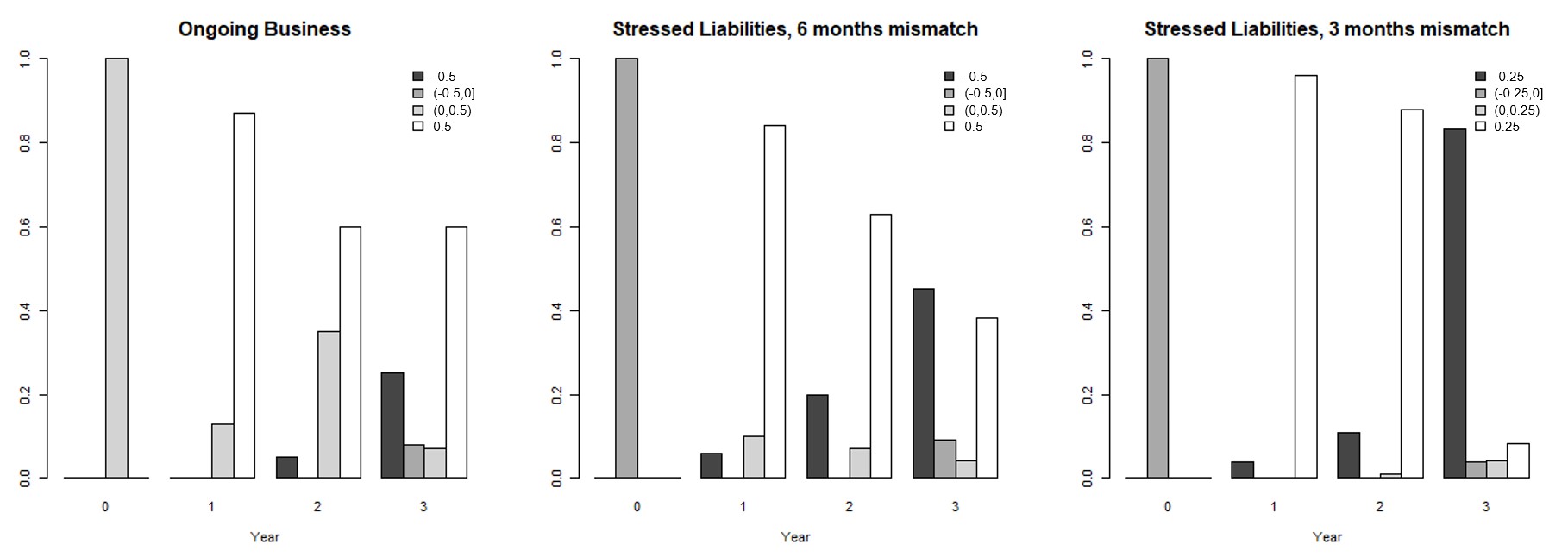}
    \caption{Percentage of nodes with specific assets-liabilities duration mismatches at different stages for the ongoing business case (left), the stressed liability case with 6-month maximum mismatch (center) and the stressed liability case with 3-month maximum mismatch (right).}
   \label{fig:Duration Comp}
\end{figure}

For duration mismatch we distinguish the following four cases, under each of the above bounds:
\begin{itemize}
    \item $-\overline{\Delta}^{(x,\lambda)}$: constraint \eqref{eq:DurationMatching} is active and liabilities duration exceed assets duration.
    \item $(-\overline{\Delta}^{(x,\lambda)};0]$: constraint \eqref{eq:DurationMatching} is not active with liabilities durations exceeding assets durations.
    \item $(0;\overline{\Delta}^{(x,\lambda)})$: constraint \eqref{eq:DurationMatching} is not active with assets durations exceeding liabilities duration. 
    \item $\overline{\Delta}^{(x,\lambda)}$: constraint \eqref{eq:DurationMatching} is active and assets durations exceed liabilities durations. 
\end{itemize}
The left panel of Figure \ref{fig:Duration Comp} shows, for the ongoing business scenario, the percentage of nodes with duration mismatching across the four groups in stages $t \in {\cal T}^{\prime}$. The evidence is of a prevalent exposure to increasing interest rates over the investment horizon.
For the stressed case with $\overline{\Delta}^{(x,\lambda)}=0.5$, we see in central panel of Figure \ref{fig:Duration Comp} that a slightly more balanced duration matching takes place at late stages.

We also consider for the stressed case a tightening of the duration matching constraint to $[-0.25;0.25]$ (3 months). 
This was motivated by observing that in the other instances most of the nodes showed a duration mismatch in the extreme values $-0.5$ and 0.5. 
As expected, by reducing to 0.25 the maximum duration mismatch, most of the nodes still lie at the boundaries of the domain, as illustrated in the right panel of Figure \ref{fig:Duration Comp}. 
Regarding portfolio performances, the tighter duration constraint strengthens the incentive to invest in fixed income assets to meet the stricter matching requirements. 
However, this comes at the cost of poorer overall performance, as illustrated in line 8 of Table \ref{Prof}.

\subsection{Multi-period stochastic dominance and funding conditions} \label{sec:SDimpact}

\begin{figure}[ht!]
    \centering
    \includegraphics[width=\textwidth]{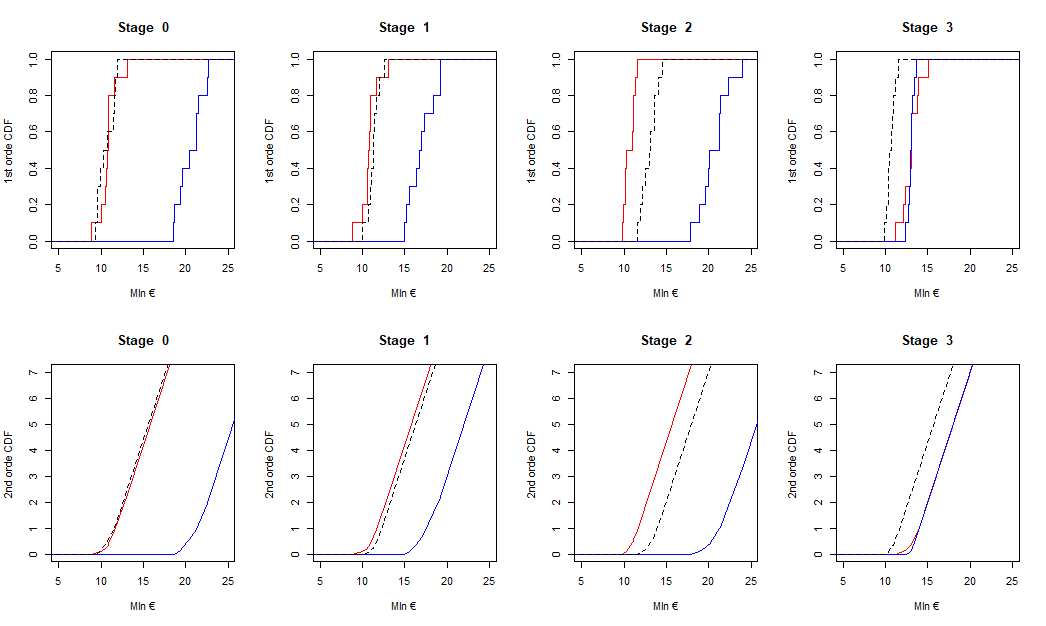}
    \caption{First order and second order CDFs for liabilities (red solid lines) and asset portfolios with (solid blue line) and without (dotted black line) SD constraints in representative nodes at different stages in the ongoing business scenario.}
    \label{fig:SD_fig1}
\end{figure}

\begin{figure}[ht!]
    \centering
    \includegraphics[width=\textwidth]{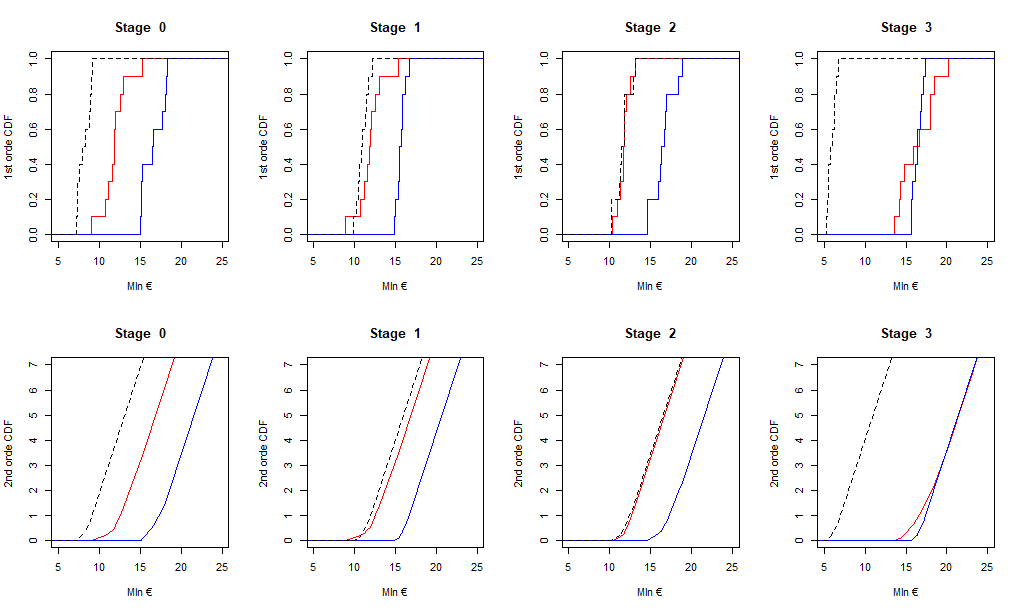}
    \caption{First order and second order CDFs for liabilities (red solid lines) and asset portfolios with (solid blue line) and without (dotted black line) SD constraints in representative nodes at different stages in the stressed liability scenario.}
    \label{fig:SD_fig2}
\end{figure}

In our numerical experiments, we test the quantile function decomposition method using the event cuts as described in Section \ref{sec:decomposition}.
According to Theorem \ref{theorem1}, the stochastic ordering is only imposed for nodes of the scenario tree at stage $T-1$. 
However, it guarantees that the asset portfolio stochastically dominates the liabilities at each stage of the problem.
Figures \ref{fig:SD_fig1} and \ref{fig:SD_fig2} provide graphical evidence for the ongoing business and the stressed liability cases, respectively.
Specifically, in these figures, we consider a representative node at each stage $t \in \mathcal{T}^\prime$ and compare first-order and second-order CDFs of three distributions in its children nodes: the distribution of liabilities (solid red line) and the asset portfolio with (solid blue line) and without (dotted black line) stochastic dominance constraints.
The representative node at stage $3$ is chosen as one where the SD constraint is actively enforced in both liability cases, making its impact most evident. 
At intermediate stages $t, \, 1\leq t \leq T-2$, the representative node is the ancestor of the corresponding node at stage $t+1$. 
At stage 0, we consider the root node.

As observed, second order stochastic dominance constraints shift the CDFs for the assets portfolio to the right, ensuring dominance over the liabilities distribution. 
This effect is particularly pronounced in the stressed liability case, where the absence of SD constraints results in solvency risks for the financial intermediary at most stages.
Additionally, Figures \ref{fig:SD_fig1} and \ref{fig:SD_fig2} confirm the validity of Theorem \ref{theorem1}, demonstrating that once stochastic dominance is enforced at stage $T-1$, it holds across all previous stages. 
When SD constraints are imposed, the model increases the initial investment $K_{0}$ in order to purchase more assets at $t=0$, raising portfolio value across all stages. 
For instance, in the ongoing business scenario, as reported in Table \ref{Root Comp}, without dominance constraints, the optimal initial investment $K_{0}$ is $9.095$ million \texteuro{}, $42.4\%$ lower than in the solution with stochastic dominance constraints, not allowing to fully cover liabilities.
The impact is even more pronounced in the stressed liability case, where the optimal investment falls from 18.660 million \texteuro{} to 7.748 million \texteuro{}.
Examining lines 9 and 10 of Tables \ref{Debt} and \ref{Prof}, we observe that excluding SD constraints reduces the initial capital requirement, which is managed more efficiently, leading to higher profits and lower borrowing needs. 
However, this approach fails to ensure stable solvency over the planning period, as the assets portfolio is unable to fully cover liabilities, leaving the financial intermediary exposed to significant solvency risks.
%
%

To further illustrate the impact of SD constraints, we analyse how varying liability funding requirements $\varphi$ affect the optimal portfolio strategy, considering under-funding and over-funding cases for both the ongoing business case and stressed liability scenarios.   
Lines 11 to 14 of Tables \ref{Debt}, \ref{Prof}, and \ref{Root Comp} summarize the results.
In both cases, as $\varphi$ increases, initial capital rises, profits decline, and debt increases, particularly in the late stages. 
Finally, to illustrate the significant impact of second order stochastic dominance constraints on the optimal investment policy, Table \ref{tab:Percentage_SD} reports the percentage of nodes at stage $T-1$ where these constraints are active across different instances. 
\begin{table}[ht!]
    \centering
    \begin{tabular}{llrr}
    \hline
        $\varphi$ & Covering Strategy & Ongoing Business  & Stressed Liability \\
        \hline
        0.8 & under-funding & 7.3\% & 10.4\% \\
         1 & full-funding & 19.3\% & 31.6\% \\
         1.1 & over-funding & 24.7\%  & 34.3\% \\
         \hline
    \end{tabular}
    \caption{Percentage of nodes at stage $T-1$ where the stochastic dominance constraint is active in the ongoing business and stressed liability case for different values of the funding ratio $\varphi$.}
    \label{tab:Percentage_SD}
\end{table}
Under a full-funding liabilities coverage, in the optimal solution of the ongoing business scenario, stochastic dominance constraints are active in $193$ nodes (i.e., 19.3\% of the total number of nodes on which these constraints are imposed). 
In the stressed liability scenario, increasing the annual liability growth rate $\mu_{1,\xi}$ from 1\% to 5\% and the volatility $\sigma_{1,\xi}$ from 3\% to 5\% leads to a higher number of active stochastic dominance constraints, reaching 31.6\% compared to 19.3\% in the base case. 
As expected, a higher funding ratio results in a greater proportion active constraints, highlighting their crucial role in maintaining a stable solvency condition for the financial intermediary.

\section{Conclusions}
\label{Conclusions}

In this paper, we provide a novel formulation for a long-term ALM problem under interest rate, inflation and credit risk exposure, with solvency and funding protection. The model here formulated jointly manages the initial capital injection, the control of debt and the internal profit generation. Thus, the solution to such a model enables the ALM managers to preserve the company funding status and fulfill liability obligations, while spanning different risk profiles. 
By enforcing in a time-consistent manner SSD of the portfolio return distribution over liabilities' distribution, our model results in a robust optimal risk-averse policy, preserving the intermediary's financial equilibrium.
The proposed model represents a significant extension of a practically and operationally relevant ALM model for a large insurer, as, for the first time in the literature, stochastic-ordering relations and dynamic risk measures are included in a multistage framework in a time-consistent manner.
Furthermore, we demonstrate that imposing the SSD constraint in a time-consistent manner at the last-but-one stage is sufficient to enforce the SSD ordering at each stage.

To tackle the problem, we develop an efficient decomposition scheme and discuss its convergence. Specifically, by using the dynamic programming formulation of the problem, we propose a new recursive solution approach based on a version of the multi-cut method in which additional cuts approximate the stochastic order constraints and the risk measures in the objective function.

The proposed methodology is validated on a case study inspired by an European insurance intermediary over a 5-year planning horizon, with portfolio rebalancing occurring in five stages, assuming either a base or a stressed liability scenario. Computational results show the effectiveness of the proposed method, which converges to the optimal solution in 7 iterations with a scenario tree with 11111 nodes and 10000 scenarios. Moreover, the stressed scenario does not affect the computational performances of the proposed method.
From a financial perspective, we notice how the stressed condition for liabilities implies a significant increase of the dedicated capital, which would make such a scenario highly undesirable due to the high cost of capital. However, the base and the stressed liability scenarios show a similar evolution for the assets portfolio, which is managed so as to limit the debt, while pursuing investment profits. 

A post-optimality analysis based on a sensitivity of the weights of debt, profits and initial invested capital shows that, when more relevance is assigned to the investment profits in the objective function, the cumulative debt is increased, especially in the early stages, and the portfolio is managed to attain higher investment profits. 
Higher penalties for the initial investment determine a significant reduction of the invested capital, requiring the availability of extra liquidity in the late stages, thus increasing the debt.
We further assess the interest rate risk exposure through the duration-matching constraint. Numerical experiments show a prevalent exposure of the financial intermediary to increasing interest rates over the investment horizon, having assets duration exceeding liabilities duration. 
When further tightening the duration constraint by reducing the maximum duration mismatch, we observe that the incentive to invest in fixed income assets increases but, overall, the search of superior portfolio performance leads to a relevant equity investment.

Finally, we assess the impact of SD constraints on the optimal solution. Results show that second order SD constraints
significantly affect the optimal investment policy, especially under the stressed liability scenario, raising from 19.3\% to 31.6\% of the nodes in which they are active under a full-funding liability coverage. 
Furthermore SD constraints imply an increased initial investment value in order to purchase more assets at the beginning of the investment horizon, thus rising the portfolio value over all stages to cover the liabilities.

\vspace{.5cm}
{\bf Acknowledgements}: Giorgio Consigli acknowledges the support from Khalifa University of Science and Technology, Grant FSU2022-010 award 000634-00001, project no. 8474000393.
\\ Darinka Dentcheva acknowledges the support of the Air Force Office of Scientific Research with the award FA9550-24-1-0284.
\\ Francesca Maggioni acknowledges the support from the PRIN2020 project "ULTRA OPTYMAL - Urban Logistics
and sustainable TRAnsportation: OPtimization under uncertainTY
and MAchine Learning", funded by the Italian
University and Research Ministry (Grant 20207C8T9M,  https://ultraoptymal.unibg.it) and from Gruppo Nazionale per il Calcolo Scientifico (GNCS-INdAM).
\\ Giovanni Micheli acknowledges the support from Gruppo Nazionale per il Calcolo Scientifico (GNCS-INdAM).

\vspace{.3cm}
{\bf Declarations}: The authors have no competing interests to declare that are relevant to the content of this article. The data that support the findings of this study are openly available at \url{http://finance.yahoo.com}, and upon request to the authors.

\bibliography{stochdom}

\newpage
\begin{appendices}
\setcounter{table}{0}
\renewcommand{\thetable}{A\arabic{table}}
\section{Statistical models parameters estimation}
\label{Stat}
In this appendix we present results for the estimation by OLS method of the coefficients associated with the statistical models introduced in Section \ref{sec:statmodel}. Input data are monthly observations from December 2018 to December 2022.
The statistical evidences for yield curve parameters are presented in Tables \ref{Est NSS} and \ref{Var}. Specifically, Table \ref{Est NSS} provides statistics for the estimated coefficients $b_{j,t}^y$, $j=1,2,3$ and $\gamma_t$ of the \textit{Nelson-Siegel-Svensson} model \eqref{eq:NSmodel}, while Table \ref{Var} shows the symmetric variance and covariance matrix of coefficients $b_{j,t}^y$, $j=1,2,3$, which is needed in order to implement the arbitrage free calibration method described in \cite{Chris2009}. 
Table \ref{Est Coeff} provides estimates and standard errors for coefficients $a_{j}^y,\ j=0,1,2,3$ of the decay factor model \eqref{eq:LambdaModel}, while estimations for coefficients  $\alpha^\pi$ and $\sigma^\pi$ of the inflation model \eqref{eq:CPI} and $c_{j},\ j=0,1,2$ of the credit spread model \eqref{eq:spread} are reported in Table \ref{Est Coeff2}.
\begin{table} [ht!]
\begin{small}
\begin{tabular}{|lrrrr|} \hline
&$b^y_{1,t}$& $b^y_{2,t}$ &$b^y_{3,t}$& $\gamma_t$\\ \hline
Mean&$0.0247$& $-0.0188$&$0.0182$&$4.9924$ \\
Std&$0.0162$& $0.0097$&$0.0111$&$2.7213$ \\
25\%-Q&$0.0111$& $-0.0254$&$0.0114$&$3.6704$ \\
50\%-Q& $0.0263$& $-0.0185$& $0.0184$&$4.5455$\\
75\%-Q&$0.0403$&  $-0.0114$&$0.0269$&$5.6551$ \\ \hline
\end{tabular}
\caption{\label{Est NSS} Mean, standard deviation (Std), first quartile (25\%-Q), median (50\%-Q), and third quartile (75\%-Q) of estimated coefficients $b_{j,t}^y,\ j=1,2,3$ and $\gamma_t$ of the yield curve model \eqref{eq:NSmodel}.}
\end{small}
\end{table}
\begin{table} [ht!]
\begin{small}
\begin{tabular}{|lrrr|} \hline
&$b^y_{1,t}$& $b^y_{2,t}$ &$b^y_{3,t}$\\ \hline
$b^y_{1,t}$&$0.000257$&$-$ &$-$ \\
$b^y_{2,t}$&$-0.000063$& $0.000093$&$-$\\
$b^y_{3,t}$&$0.000012$&  $-0.000062$&0.000144\\ \hline
\end{tabular}
\caption{\label{Var}Variance and covariance symmetric matrix for coefficients $b_{j,t}^y, \ j=1,2,3$.}
\end{small}
\end{table}
\begin{table}
\begin{small}
\begin{tabular}{|l|rrrr|}
\hline
 & \multicolumn{4}{c|}{Decay factor \eqref{eq:LambdaModel}} \\ 
&$a^y_{0}$& $a^y_{1}$&$a^y_{2}$& $a^y_3$ \\ \hline
Mean&$7.0549$& $47.7621$&$121.3425$&$50.8006$ \\
Std &$0.3831$&  $11.0799$&$21.8113$&$16.4272$ \\ \hline
\end{tabular}
\caption{\label{Est Coeff}
{Estimates and standard errors (Std) for coefficients $a_{j}^y,\ j=0,1,2,3$ of the decay factor model \eqref{eq:LambdaModel}. All coefficients are statistically significant at
the 1\% level.}}
\end{small}
\end{table}
\begin{table}
\begin{small}
\begin{tabular}{|l|rr|rrr|}
\hline
 & \multicolumn{2}{c|}{Inflation \eqref{eq:CPI}}& \multicolumn{3}{c|}{Credit Spread \eqref{eq:spread}}\\ 
& $\alpha^\pi$& $\sigma^\pi$
& $c_{0}$& $c_{1}$&$c_2$
\\ \hline
Mean& $0.2344$*& $0.0508$*
& $0.0614$& $0.9479$*&$4.1689$*
\\
Std& $0.0172$& $0.0020$& $0.0351$& $0.0189$&$1.1905$\\ \hline
\end{tabular}
\caption{\label{Est Coeff2} Estimates and standard errors (Std) for coefficients $\alpha^\pi$ and $\sigma^\pi$ of the inflation model \eqref{eq:CPI}, and $c_{j},\ j=0,1,2$ of the credit spread model \eqref{eq:spread}. Statistical significance at
the 1\% level is denoted by ’*’.}
\end{small}
\end{table}

With regard to asset price returns, Table \ref{Est Asset} provides estimates of the regression coefficients $b_i=(b_{i,0}, b_{i,1}, b_{i,2}, b_{i,3}, b_{i,4})^\top, \ i \in {\cal I}$, in models \eqref{eq:A1model}, \eqref{eq:A2model}, \eqref{eq:A3model} and \eqref{eq:A4model}, and the corresponding coefficients of determination $R^2$.
\begin{table}[ht!]
\begin{small}
\begin{tabular}{|lrrrrrr|}
\hline
Asset & $b_{i,0}$  & $b_{i,1}$ & $b_{i,2}$ & $b_{i,3}$ & $b_{i,4}$
& $R^2$ \\
\hline
\multirow{2}{*}{Smart.MI} & $0.0001$ & $0.2997$* & $0.0398$* & $-0.0032$ & & \multirow{2}{*}{$0.82$} \\
 & (0.0001) & (0.1096) & (0.0136) & (0.038) & & \\
\hline
\multirow{2}{*}{IBGS.L} & $0.0324$* & $-0.1358$ & $4.4568$ & $-0.4895$ & & \multirow{2}{*}{$0.79$} \\
 & (0.0118) & (0.1154) & (0.7706) & (0.0607) & & \\
\hline
\multirow{2}{*}{DBXR.DE} & $0.0029$* & $-0.0354$ & $0.4042$ & $-0.2334$* & & \multirow{2}{*}{$0.74$} \\
 & (0.0005) & (0.1169) & (0.3005) & (0.0212) & & \\
\hline
\multirow{2}{*}{SYBV.DE} & $-0.0007$ & $0.0048$ & $0.3902$* & $0.3065$* & & \multirow{2}{*}{$0.68$} \\
 & (0.0038) & (0.1322) & (0.0614) & (0.0842) & & \\
\hline
\multirow{2}{*}{LQD} & $0.0019$ & $-0.0156$ & $-0.1251$* & $-0.0757$ & & \multirow{2}{*}{$0.59$} \\
 & (0.0023) & (0.1186) & (0.0468) & (0.5679) & & \\
\hline
\multirow{2}{*}{BND} & 0.0006 & $-0.0113$ & 0.0943 & $-0.9348$* & & \multirow{2}{*}{$0.74$} \\
 & (0.0011) & (0.0701) & (0.0537) & (0.1856) & & \\
\hline
\multirow{2}{*}{TDTT} & $-0.0012$ & 0.0573 & $-0.0146$ & 0.1298* & 0.0302* & \multirow{2}{*}{$0.68$} \\
 & (0.0044) & (0.1119) & (0.4360) & (0.0470) & (0.0112) & \\
\hline
\multirow{2}{*}{TDTF} & $-0.0021$ & 0.0479 & $-0.1851$* & 0.1544 & 0.0292* & \multirow{2}{*}{$0.71$} \\
 & (0.0032) & (0.1125) & (0.0847) & (0.2089) & (0.01002) & \\
\hline
\multirow{2}{*}{IMEU.AS} & $-0.0380$ & $-0.1206$ & $-7.0718$* & $0.9165$* & $-0.8562$ & \multirow{2}{*}{$0.66$} \\
 & (0.0316) & (0.1156) & (2.3629) & (0.4133) & (1.3166) & \\
\hline
\multirow{2}{*}{IWM} & $-0.0089$ & $-0.0422$ & $-5.4131$* & $-1.1704$ & $-2.0308$* & \multirow{2}{*}{$0.64$} \\
 & (0.0456) & (0.1178) & (3.2857) & (0.9410) & (0.7056) & \\
\hline
\multirow{2}{*}{EEM} & $-0.0252$ & $0.0154$ & $-5.7660$* & $-1.2073$ & $0.8823$ & \multirow{2}{*}{$0.63$} \\
 & (0.0348) & (0.1156) & (2.2189) & (1.4678) & (1.4537) & \\
\hline
\multirow{2}{*}{EWJ} & 0.0044 & 0.0590 & $-0.2506$ & $-0.6725$ & 0.0341 & \multirow{2}{*}{$0.36$} \\
 & (0.0060) & (0.0688) & (0.2014) & (0.7057) & (0.3794) & \\
\hline
\multirow{2}{*}{SUA.AS} & $-0.0141$ & $-0.1397$* & $-0.6447$ & 0.9519 & 0.4534 & \multirow{2}{*}{$0.59$} \\
 & (0.0475) & (0.0495) & (1.1904) & (1.2369) & (5.0028) & \\
\hline
\multirow{2}{*}{EURUSD} & $-0.0012$ & $-0.0072$ & 0.0538 & 0.1969* & & \multirow{2}{*}{$0.25$} \\
 & (0.0021) & (0.0685) & (0.1211) & (0.0747) & & \\
\hline
\end{tabular}
\caption{\label{Est Asset} Estimates and in parenthesis standard errors of regression coefficients $b_i=(b_{i,0}, b_{i,1}, b_{i,2}, b_{i,3}, b_{i,4})^\top$,  $i \in {\cal I}$. Statistical significance at the 1\% level is denoted by '*'.}
\end{small}
\end{table}

\section{Scenario generation algorithm}
\label{App_Scen}
\setcounter{figure}{0}
\renewcommand{\thefigure}{B\arabic{figure}}
\setcounter{table}{0}
\renewcommand{\thetable}{B\arabic{table}}

Let $\theta$ be the vector including all statistical coefficients estimated by OLS method according to the models from \eqref{eq:NSmodel} to \eqref{eq:gaincoeff} and all the parameters of the yield curve model, the inflation and the credit spread process specified at $t=0$. 

Let $\xi_n:=\{r_{i,n},g_{i,n},\lambda_{j,n},L_{j,n},c_n,\delta_{j,n}^\lambda\}$ be a coefficient tree process on the node $n \in {\cal N}$ of the scenario tree. This vector will then include all the random parameters specified in the ALM model.
Since returns $r_{i,n}$ of assets $i \in {\cal I}$ in node $n$ and values $\lambda_{j,n}$ of liability $j \in {\cal J}$ in node $n$ depend on the yield curve $y_{t,\tau}$, the inflation process $\pi_t$ and the credit spread process $s_t^{IG}$, values for the vector $\xi_n$ are generated by applying a two-step procedure, with the first step being the generation of a random vector process $\omega_n:=(y_{n,\tau},\pi_n,s_n^{IG})$ for yield curve rates, inflation and credit spread, and the second step being the generation of the stochastic ALM model coefficients.
Values for the random vector process $\omega_n$, referred to as the \textit{core economic model}, are determined by applying Algorithm \ref{algo:scengen1}. Specifically, the input to the algorithm is represented by the vector $\theta$. According to the planning horizon $T$, to the stage composition ${\cal T}$ and to the branching degree vector, the \textit{Nodal Partition Matrix} (NPM) is generated. Such a matrix has $N_T$ rows and $T+1$ columns. Each row of the matrix is a scenario for the core economic model, determined by applying models \eqref{eq:NSmodel}, \eqref{eq:CPI} and \eqref{eq:spread} with monthly increments from $t=0$ to $T$.
\begin{algorithm} 
\begin{small}
\begin{flushleft}
\textbf{Input} Vector $\theta$: \\ \hspace{1cm}(a) Nelson-Siegel-Svensson parameters $b^y_{1,0}, b^y_{2,0}, b^y_{3,0}, \lambda_0$ plus stochastic coefficients for $b^y_{j,t}$, $j=1,2,3$ and for $\lambda_t=\lambda(b^y_{j,t})$.\\ \hspace{1cm}(b) inflation process coefficients $\alpha^{\pi}, \sigma^{\pi}$ and $\pi_0$.\\ \hspace{1cm}(c) credit spread process coefficients $c_0, c_1, c_2$ and initial condition $s^{IG}_0$.\\ \hspace{1cm}(d) errors distributions for each model.
\begin{enumerate}
\item Specify planning horizon $T$ and stage composition ${\cal T}$. 
\item Generate the \textit{Nodal Partition Matrix} (NPM) of ${\cal N}_T$ rows and $T+1$ columns. 
\item \textbf{For} $t=1:T$
\begin{itemize}
\item[] {\bf For} $n \in {\cal N}_t$ 
\begin{itemize}
    \item[] {\bf For} $h=(t_{a(n)},t_{a(n)}+\Delta t,...,t_{n}-\Delta t, t_n)$,    \textit{$\Delta t$ monthly increments between nodes}
    \begin{itemize}
    \item[-] generate yield curve inter-stage increments from \eqref{eq:NSmodel}:  $y_{h,\tau}=y(b^y_{1,h},b^y_{2,h},b^y_{3,h},\lambda_h)$;
    \item[-] generate inflation increments \eqref{eq:CPI} $\pi_h=\pi(\alpha^{\pi},\sigma^{\pi})$;
    \item[-] generate credit spread increments from \eqref{eq:spread} $s^{IG}_h=s(c_0,c_1,c_2)$.
    \end{itemize}
    \item[]  {\bf End For}
    \end{itemize}
    \item[] {\bf End For}
    \end{itemize}
\item[] \textbf{End For}
\end{enumerate}
\textbf{Output} $y_{n,\tau}$, $\pi_{n}$, $s_n^{IG}$ scenario paths.
\end{flushleft}
 \caption{Scenario generation - core economic model}
 \label{algo:scengen1}
 \end{small}
\end{algorithm}

The scenarios of the core economic model are input to the asset returns and liability costs models, which determine the scenarios of the coefficient process $\xi_n$, as detailed in Algorithm \ref{algo:scengen2}.
The initial conditions of the scenario generation are defined at the root node by $\xi_0:=\{r_{i,0},g_{i,0},\lambda_{j,0},L_{j,0},c_0\}$. We distinguish here between the investment horizon $T$ and the liability valuation horizon $T_{\lambda}$: this term reflects the number of years in the future in which liabilities are accounted for. 
The algorithm to determine values for $\xi_n$ consists of a forward pass and a backward pass. Similarly to the previous algorithm, in the forward pass we generate scenarios for parameters $r_{i,n}$, $g_{i,n}$, $L_{j,n}$, and $c_{n}$ by applying the statistical models introduced in Sections \ref{sec:liabdyn} and \ref{sec:assetdyn} with monthly increments from $t=0$ to $T+T_\lambda$. We then determine with a backward recursion the stochastic values of parameters $\lambda_{j,n}$ and $\delta_{j,n}^\lambda$ by discounting the expected payments according to equations \eqref{eq:lambda1} and \eqref{eq:duration_liab}. Finally, the total liability value in each node $\Lambda_n=\sum_{j \in {\cal J}} \lambda_{j,n}$ is computed.

\begin{algorithm} [ht!]
    \begin{small}
        \begin{flushleft}
      \textbf{Input} Initial conditions for $\xi_0:=\{r_{i,0}=g_{i,0}=0, \lambda_{j,0},L_{j,0},c_0\}$. \\ 
      \textbf{Input} ALM horizon $T$, stage composition ${\cal T}$ and liability evaluation horizon $T_{\lambda}$. \\ 
     \textbf{Input} Initial term structure of interest rates $y_{0,\tau}$.  \\
     \begin{enumerate}
      \item Generate NPM for liability evaluation: ${\cal N}_T$ rows and $T + {T}_{\lambda}+1$ columns.
      \item {\bf For} $t=1:(T + T_{\lambda})$ \ \ \textit{forward pass}
      \begin{itemize}
      \item[]  {\bf For} $n \in {\cal N}_t$   \\
      \begin{itemize}
          \item[] {\bf For} $h=(t_{a(n)},t_{a(n)}+\Delta t,...,t_{n}-\Delta t, t_n)$ \\
          \begin{itemize}
                \item[] Compute $\{r_{i,h}, g_{i,h}, L_{j,h}, c_h\}$  \ \ from   \eqref{eq:netrevenues}, \eqref{eq:rev_inc}, \eqref{eq:netpay}, \eqref{eq:pay_inc}, \eqref{eq:A1model}, \eqref{eq:A2model}, \eqref{eq:A3model}, \eqref{eq:gaincoeff}
        \end{itemize}
                \item[] {\bf End For}
      \end{itemize}
      \item[] {\bf End For}
      \end{itemize}
      \item[] {\bf End For}
      \item {\bf For} $t=(T + T_{\lambda}) : 0$ \  \  \  \  \textit{backward recursion}  \\
      \begin{itemize}
      \item[] {\bf For} $n \in {\cal N}_t$  \\
      \begin{itemize}
      \item {\bf For} $j=1,2,..,J$ \\
      \begin{itemize}
      Compute $\lambda_{j,n}$ and $\delta_{j,n}^\lambda$ \ \ from \eqref{eq:lambda1} and \eqref{eq:duration_liab}.
      \end{itemize}
      \item[] {\bf End For} \\ 
      \item Compute $\Lambda_n=\sum_{j \in {\cal J}} \lambda_{j,n}$.
      \end{itemize}
      \item[] {\bf End For}
      \end{itemize}
      \item[] {\bf End For}
     \end{enumerate}
     \textbf{Output} $r_{i,n}, g_{i,n}, L_{j,n}, c_n, \lambda_{j,n}, \delta_{j,n}^\lambda,\Lambda_n$    
        \end{flushleft}
        \caption{Scenario generation - coefficient process}
        \label{algo:scengen2}
    \end{small}
\end{algorithm}

\begin{figure}[ht!]
    \centering
    \includegraphics[width=1\linewidth]{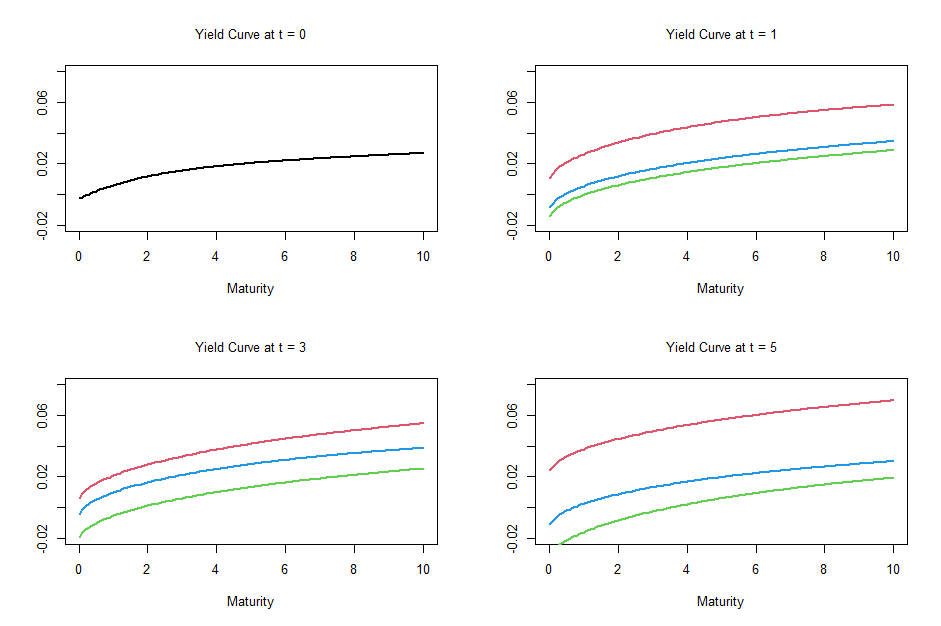}
    \caption{Graphical representation of the yield curve at stage $t=0$ and the evolution of yield curves in selected scenarios at years 1, 3, and 5.}
    \label{fig:YC}
\end{figure}

 \vspace{.2cm} 
In the following, we present the simulation outputs of the scenario generation algorithms. 
Specifically, Table \ref{tab:stats_NS} compares statistical measures computed on historical and simulated data for the \textit{Nelson-Siegel-Svensson} parameters \(b_{j,t}^y\), with \(j=1,2,3\).  
\begin{table}[ht!]
    \centering
    \begin{tabular}{l|l|rrrrrr}
    \hline
 Parameter & Series & Mean& Std& Skewness& Kurtosis& 5\%-Q& 1\%-Q
\\
\hline
  \multirow{2}{*}{$b_{1,t}^y$} &History& 0.0247& 0.0162& --0.0471& --1.4243
& 0.0005
& --0.0012
\\
  &Simulation& 0.0183& 0.0052& 0.1400& --0.6436
& 0.0033
& 0.0021
\\
\hline
  \multirow{2}{*}{$b_{2,t}^y$}  &History&  --0.0188&  0.0097&  --0.4426&  --0.5550
&  --0.0387
&  --0.0410
\\
          &Simulation&  --0.0111&  0.0015&  --0.1239&  --0.5432
&  --0.0225
&  --0.0243
\\
\hline
 \multirow{2}{*}{$b_{3,t}^y$} &History&  0.0182&  0.0117&  0.2984&  --0.5502
&  0.0000
&  0.0000
\\
          &Simulation&  0.0111&  0.0177&  0.1114&  --0.5317
&  0.0018
&  0.0010
\\
\hline
    \end{tabular}
     \caption{Mean, standard deviation (Std), skewness, kurtosis, 5\%-quantile (5\%-Q), 1\%-quantile (1\%-Q) of \textit{Nelson-Siegel-Svensson} parameters $b_{j,t}^y, \, j=1,2,3$: historical versus simulation evidence at the end of year 5.}
    \label{tab:stats_NS}
\end{table}
Since the generation of reliable yield curves is a key requirement for an ALM model, Figure \ref{fig:YC} illustrates the evolution of simulated yield curves. 
To ensure consistency over time, we compare the yield curve at stage \(t = 0\) with those at years 1, 3, and 5 under three representative scenarios: (i) the scenario with the lowest 3-month interest rate (green line), (ii) the scenario with the median 3-month interest rate (blue line), and (iii) the scenario with the highest 3-month interest rate (red line).

\begin{figure}[ht!]
    \centering
    \includegraphics[width=0.9\linewidth]{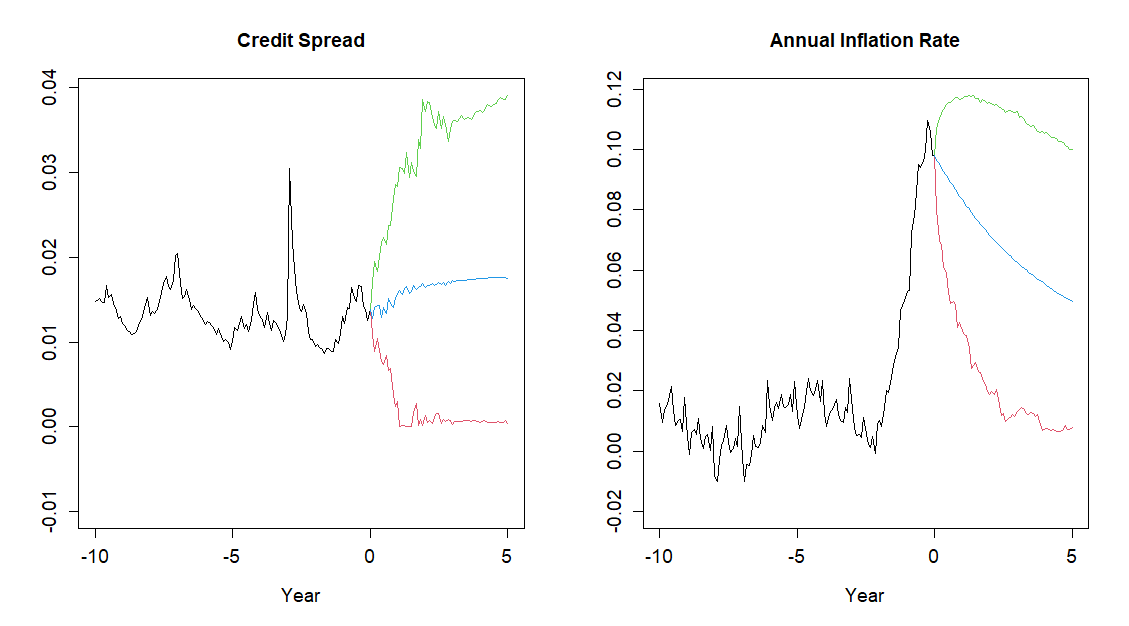}
    \caption{Evolution of credit spread (left) and inflation rate (right) over time: historical versus simulation evidence.}
    \label{fig:Spread + Infl}
\end{figure}

Figure \ref{fig:Spread + Infl} presents the historical 10-year data and 5-year simulation outputs for the investment-grade credit spread (\(s_t^{IG}\)) and the inflation rate (\(\pi_t\)), with summary statistics reported in Table \ref{tab:stats_infl}. 
To provide a comprehensive view of the simulated variability, Figure \ref{fig:Spread + Infl} highlights, for each projected time step, the minimum, maximum, and median values of the generated scenarios.  
\begin{table}
    \centering
    \begin{tabular}{l|l|rrrrrr}
    \hline
        Parameter &  Series & Mean &  Std&  Skewness&  Kurtosis&  5\%-Q&  1\%-Q
\\
\hline 
         \multirow{2}{*}{$s_t^{IG}$} & History&  0.0153&  0.0085&  3.1484&  15.9086&  0.0078&  0.0057
\\
         & Simulation&  0.0168&  0.0089&  0.3014&  2.8803&  0.0066&  0.0015
         \\
\hline
         \multirow{2}{*}{$\pi_t$} & History&  0.0333&  0.0339&  1.0568&  2.5990&  0.0022&  0.0007
\\
         & Simulation&  0.0526&  0.0330&  0.5614&  3.3921&  0.0281&  0.0194
\\
\hline
    \end{tabular}
    \caption{Mean, standard deviation (Std), skewness, kurtosis, 5\%-quantile (5\%-Q), 1\%-quantile (1\%-Q) of credit spread $s_t^{IG}$ and annual inflation rate $\pi_t$: historical versus simulation evidence at the end of year 5.}
    \label{tab:stats_infl}
\end{table}

\begin{figure}[ht!]
    \centering
    \includegraphics[width=1\linewidth]{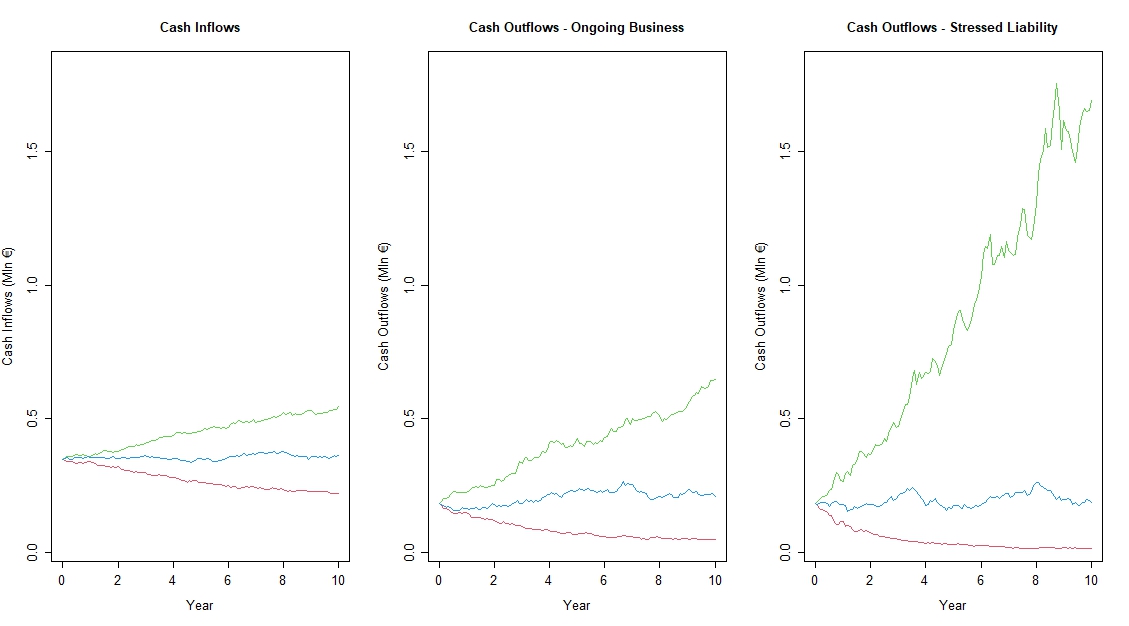}
    \caption{10-year simulation outputs for cash outflows $L_{1,t}$ and cash inflows $c_t$.}
    \label{fig:Cash IN-OUT}
\end{figure}

Finally, scenarios for cash inflows and outflows under both ongoing business and stressed liability conditions are depicted in Figure \ref{fig:Cash IN-OUT}, while the corresponding statistics are summarized in Table \ref{tab:liab_stats}. 
In this case, the projection horizon extends to \(T + T_{\lambda} = 10\) years to ensure a proper evaluation of liabilities. 
Similar to previous analyses, we highlight the total variability of simulated data by showing the minimum, maximum, and median values at each time step.  
The correspondent statistics are reported in Table \ref{tab:liab_stats}.
\begin{table}[ht!]
    \centering
    \begin{tabular}{l|rrrrrr}
    \hline
 Parameter & Mean& Std& Skewness& Kurtosis& 5\%-Q& 1\%-Q
\\
\hline
Cash inflows & 0.3502 & 0.0405 & 0.3226 & 0.1345 & 0.2875 & 0.2669 \\
Cash outflows (ongoing) & 0.1847 & 0.0666 & 1.0655 & 1.9433 & 0.0973 & 0.0784 \\
Cash outflows (stressed) & 0.1957 & 0.1422 & 2.2379 & 9.4149 & 0.0513 & 0.0324 \\
\hline
    \end{tabular}
     \caption{Mean (million of euros), standard deviation (Std), skewness, kurtosis, 5\%-quantile (5\%-Q), 1\%-quantile (1\%-Q) of cash inflows and cash outflows in the ongoing business and stressed liability cases: simulation evidence at the end of year 10.}
    \label{tab:liab_stats}
\end{table}
%
%



\end{appendices}


\end{document}